\documentclass[11pt]{amsart}
\usepackage{amscd,amssymb}
\input xy
\xyoption{all}
\newcommand{\bdes}{\begin{description}}
\newcommand{\edes}{\end{description}}
\newcommand{\beqn}{\begin{eqnarray*}}
\newcommand{\eeqn}{\end{eqnarray*}}

\newcommand{\EE }{{\mathbb E}}
\newcommand{\PP }{{\mathbb P}}
\newcommand{\QQ }{{\mathbb Q}}
\newcommand{\CC }{{\mathbb C}}
\newcommand{\ZZ }{{\mathbb Z}}

\newtheorem{definition}{Definition}[subsection]
\newtheorem{theorem}{Theorem}[subsection]

\newtheorem{lemma}{Lemma}[subsection]
\newtheorem{corollary}{Corollary}[subsection]

\newtheorem{def/th}{Definition/Theorem}[subsection]
\newtheorem{remark}{Remark}[subsection]

\begin{document}
\title{Mirror symmetry for concavex vector bundles on projective spaces}
\thanks{1991 Mathematics Subject Classification. Primary 14N35. Secondary 14L30.}
\author{Artur Elezi}
\begin{abstract}
Let $X\subset Y$ be smooth, projective manifolds. Assume that
$\iota:X\hookrightarrow \PP^s$ is the zero locus of a generic
section of $V^+=\oplus_{i\in I}\mathcal O(k_i)$ where all the
$k_i$'s are positive. Assume furthermore that $\mathcal
N_{X/Y}=\iota^*(V^-)$ where $V^-=\oplus_{j\in J}\mathcal O(-l_j)$
and all the $l_j$'s are negative. We show that under appropriate
restrictions, the generalized Gromov-Witten invariants of $X$
inherited from $Y$ can be calculated via a modified Gromov-Witten
theory on $\PP^s$. This leads to local mirror symmetry on the
A-side.
\end{abstract}

\maketitle

\section{Introduction}

Let $V^+=\oplus_{i\in I}\mathcal O(k_i)$ and $V^-=\oplus_{j\in
J}\mathcal O(-l_j)$ be vector bundles on $\PP^s$ with $k_i$ and
$l_j$ positive integers. Suppose
$X\stackrel{\iota}{\hookrightarrow}\PP^s$ is the zero locus of a
generic section of $V^+$ and $Y$ is a projective manifold such
that $X\stackrel{j}{\hookrightarrow}Y$ with normal bundle
$\mathcal N_{X/Y}=\iota^*(V^-)$. The relations between
Gromov-Witten theories of $X$ and $Y$ are studied here by means of
a suitably defined equivariant Gromov-Witten theory in $\PP^s$. We
apply mirror symmetry to the latter to evaluate the gravitational
descendants of $Y$ supported in $X$.

Section $2$ is a collection of definitions and techniques that
will be used throughout this paper. In section $3$, using an idea
from Kontsevich, we introduce a modified equivariant Gromov-Witten
theory in $\PP^s$ corresponding to $V=V^+\oplus V^-$. The
corresponding $\mathcal D$-module structure (\cite
{[6]},\cite{[15]},\cite{[29]}) is computed in section $4$. It is
generated by a single function $\tilde{J}_V$. In general, the
equivariant quantum product does not have a nonequivariant limit.
It is shown in Lemma $4.1.1$ that the generator $\tilde{J}_V$ does
have a limit $J_V$ which takes values in $H^*\PP^m[[q,t]]$. It is
this limit that plays a crucial role in this work.

Let $Y$ be a smooth, projective manifold. The generator $J_Y$ of
the pure $\mathcal D$-module structure of $Y$ encodes one-pointed
gravitational descendents of $Y$. It takes values in the
completion of $H^*Y$ along the semigroup (Mori cone) of the
rational curves of $Y$. The pullback map $j^*:H^*Y\rightarrow
H^*X$ extends to a map between the respective completions. In
Theorem~\ref{thm: localthm} we describe one aspect of the relation
between pure Gromov-Witten theory of
$X\stackrel{j}{\hookrightarrow}Y$ and the modified Gromov-Witten
theory of $\PP^s$. Under natural restrictions, the pull back
$j^*(J_Y)$  pushes forward to $J_V$. It follows that although
defined on $\PP^s$, $J_V$ encodes the gravitational descendants of
$Y$ supported in $X$, hence the contribution of $X$ to the
Gromov-Witten invariants of $Y$.

The only way that $X$ remembers the ambient variety $Y$ in this
context is by the normal bundle. $Y$ can therefore be substituted
by a local manifold. This suggests that there should be a local
version of mirror symmetry (see the Remark at the end of section
$4$). This was first realized by Katz, Klemm, and Vafa
\cite{[20]}. The principle of local mirror symmetry in general has
yet to be understood. Some interesting calculations that
contribute toward this goal can be found in \cite{[8]}.

In section $5$ we give a proof of the Mirror Theorem which allows
us to compute $J_V$. A hypergeometric series $I_V$ that
corresponds to the total space of $V$ is defined. The Mirror
Theorem~\ref{thm: mthm} states that $I_V=J_V$ up to a change of
variables. Hence, the gravitational descendants of $Y$ supported
on $X$ can be computed in $\PP^s$.

Two examples of local Calabi-Yau threefolds are considered in
section $6$. For $X=\PP^1$ and $V=\mathcal O(-1)\oplus \mathcal
O(-1)$, we obtain the Aspinwall-Morrison formula for multiple
covers. If $X=\PP^2$ and $V=\mathcal O(-3)$, the quantum product
of $Y$ pulls back to the modified quantum product in $\PP^2$. The
mirror theorem in this case yields the virtual number of plane
curves on a Calabi-Yau threefold.

The rich history of mirror symmetry started in 1990 with a
surprising conjecture by Candelas, de la Ossa, Green and Parkes
(\cite{[7]}) which predicts the number $n_d$ of degree $d$
rational curves on a quintic threefold. In \cite{[15]}, Givental
presented a clever argument which, as shown later by Bini et al.
in \cite {[6]} and Pandharipande in \cite{[29]}, yields a proof of
the mirror conjecture for Fano and Calabi-Yau (convex) complete
intersections in projective spaces. Meanwhile, in a very well
written paper \cite{[26]}, Lian, Liu and Yau used a different
approach to obtain a complete proof of mirror theorem for concavex
complete intersections on projective spaces. An alternative  proof
of the convex Mirror Theorem has been given by Bertram
(\cite{[5]}). In this paper we use Givental's approach to study
the local nature of mirror symmetry and to present a proof of the
concavex Mirror Theorem.

{\bf Acknowledgements.} This work is part of author's Ph.D. thesis
at Oklahoma State University. The author would like to thank
Sheldon Katz for his passionate and tireless work in advising with
this project. Special thanks to Bumsig Kim, Carel Faber,
Ionut-Ciocan Fontanine and Zhenbo Qin who were very helpful
throughout this work. At various times the author has benefited
from conversations with Tom Graber, Rahul Pandharipande and Ravi
Vakil, to whom the author is very grateful. We would like to thank
also the referees whose help in improving this manuscript was
invaluable.
\section{Stable maps and localization}
\subsection {\bf Genus zero stable maps.} Let $\overline
M_{0,n}(X,\beta)$ be the Deligne-Mumford moduli stack of pointed
stable maps to $X$. For an excellent reference on the construction
and its properties we refer the reader to \cite {[13]}. We recall
some of the features on $\overline M_{0,n}(X,\beta)$ and establish
some notation. For each marking point $x_i$ let $e_i:{\overline
M}_{0,n}(X,\beta)\rightarrow X$ be the evaluation map at $x_i$ and
${\mathcal L}_i$ the cotangent line bundle at $x_i$. The fiber of
this line bundle over a moduli point $(C,x_1,...,x_n,f)$ is the
cotangent space of the curve $C$ at $x_i$. Let $\pi_k:
\overline{M}_{0,n}(X,\beta)\rightarrow {\overline
M}_{0,n-1}(X,\beta)$ be the morphism that forgets the $k$-th
marked point. The obstruction theory of the moduli stack
${\overline M}_{0,n}(X,\beta)$ is described locally by the
following exact sequence
\[0\rightarrow \text{Ext}^0(\Omega_C(\sum_{i=1}^{n}x_i),\mathcal O_C)\rightarrow \text{H}^0(C,f^*TX)\rightarrow \mathcal T_M\rightarrow\]
\begin{equation}
 \rightarrow \text{Ext}^1(\Omega_C(\sum_{i=1}^{n}x_i),\mathcal O_C)\rightarrow \text{H}^1(C,f^*TX)\rightarrow \Upsilon\rightarrow 0. \label {seq: deformation}
\end{equation}
(Here and thereafter we are naming sheaves after their fibres). To
understand the geometry behind this exact sequence we note that
$\mathcal T_M=\text{Ext}^1(f^*\Omega_X\rightarrow
\Omega_C,\mathcal O_C)$ and
$\Upsilon=\text{Ext}^2(f^*\Omega_X\rightarrow \Omega_C,\mathcal
O_C)$ are respectively the tangent space and  the obstruction
space at the moduli point $(C,x_1,...,x_n,f)$.  The spaces
$\text{Ext}^0(\Omega_C(\sum_{i=1}^{n}x_i),\mathcal O_C)$ and
$\text{Ext}^1(\Omega_C(\sum_{i=1}^{n}x_i),\mathcal O_C)$ describe
respectively the infinitesimal automorphisms and infinitesimal
deformations of the marked source curve. It follows that the
expected dimension of $\overline{M}_{0,n}(X,\beta)$ is $-K_X\cdot
\beta +\text{dim}X+n-3$.

A smooth projective manifold $X$ is called {\it convex} if
$H^1(\PP^1,f^*TX)=0$ for any morphism $f:\PP^1\rightarrow X$. For
a convex $X$,  the obstruction bundle $\Upsilon$ vanishes and the
moduli stack is unobstructed and of the expected dimension.
Examples of convex varieties are homogeneous spaces $G/P$.

In general this moduli stack may behave badly and have components
of larger dimensions. In this case, a Chow homology class of the
expected dimension has been constructed \cite{[3]} \cite{[27]}. It
is called the virtual fundamental class and denoted by
$[\overline{M}_{0,n}(X,\beta)]^{\text{virt}}$. Although its
construction is quite involved, we will be using mainly two
relatively easy properties. The virtual fundamental class is
preserved when pulled back by the forgetful map $\pi_n$. A proof
of this fact can be found in section $7.1.5$ of \cite {[9]}. If
the obstruction sheaf $\Upsilon$ is free, the virtual fundamental
class refines the top Chern class of $\Upsilon$. This fact is
proven in Proposition $5.6$ of \cite{[3]}.

\subsection{\bf Equivariant cohomology and localization theorem} The notion of
equivariant cohomology and the localization theorem is valid for
any compact connected Lie group. For a detailed exposition on this
subject we suggest Chapter $9$ of \cite{[9]}. Below we state
without proof the results that will be used in this work.

The complex torus $T=(\CC^*)^{s+1}$ is classified by the principal $T$-bundle
\begin{equation}
ET=(\CC^{\infty+1}-\{0\})^{s+1}\rightarrow BT=(\CC\PP^{\infty})^{s+1}.
\end{equation}
Let $\lambda_i=c_1(\pi_i^*(\mathcal O(1)))$ and
$\lambda:=(\lambda_0,...,\lambda_s)$. We will use $\mathcal
O(\lambda_i)$ for the line bundle $\pi_i^*(\mathcal O(1))$.
Clearly $H^*(BT)=\CC[\lambda]$. If $T$ acts on a variety $X$, we
let $X_T:=X\times_T ET.$
\begin{definition}
The equivariant cohomology of $X$ is
\begin{equation}
H^*_T(X):=H^*(X_T).
\end{equation}
\end{definition}
If $X=x$ is a point then $X_T=BT$ and $H^*_T(x)=\CC[\lambda]$. For
an arbitrary $X$, the equivariant cohomology $H^*_T(X)$ is a
$\CC[\lambda]$-module via the equivariant morphism $X\rightarrow
x$.

Let $\mathcal U$ be a vector bundle over $X$. If the action of $T$
on $X$ can be lifted to an action on $\mathcal U$ which is linear
on the fibers, $\mathcal U$ is an equivariant vector bundle and
$\mathcal U_T$ is a vector bundle over $X_T$. The equivariant
chern classes of $E$ are $c^T_k(\mathcal U):=c_k(\mathcal U_T)$.
We will use $\text{E}(\mathcal U)$ ($\text{E}_T(\mathcal U)$) to
denote the nonequivariant (equivariant) top chern class of
$\mathcal U$.

Let $X^T=\cup_{j\in J} X_j$ be the decomposition of the fixed
point locus into its connected components. $X_j$ is smooth for all
$j$ and the normal bundle $ N_j$ of $X_j$ in $X$ is equivariant.
Let $i_j:X_j\rightarrow X$ be the inclusion. The following
corollary of the localization theorem will be used extensively
here:
\begin{theorem}
Let $\alpha\in H^*_T(X)\otimes \CC(\lambda)$. Then
\begin{equation}
\int_{X_T}\alpha=\sum_{j\in
J}\int_{(X_j)_T}\frac{i_j^*(\alpha)}{\text{E}_T(N_j)}.
\end{equation}
\end{theorem}

A  basis for the characters of the torus is given by
$\varepsilon_i(t_0,...,t_s)=t_i$. There is an isomorphism between
the character group of the torus and $H^2(BT)$ sending
$\varepsilon_i$ to $\lambda_i$. We will say that {\it the weight}
of the character $\varepsilon_i$ is $\lambda_i$.

For an equivariant vector bundle $\mathcal U$ over $X$ it may
happen that the restriction of $\mathcal U$ on a fixed point
component $X_j$ is trivial (for example if $X_j$ is an isolated
point). In that case $\mathcal U$ decomposes as a direct sum
$\oplus_{i=1}^{m}\mu_i$ of characters of the torus. If the weight
of $\mu_i$ is $\rho_i$, then the restriction of $c^T_k(\mathcal
U)$ on $X_j$ is the symmetric polynomial
$\sigma_k(\rho_1,...,\rho_m)$.

Our interest here is for $X=\PP^s$. For any action of $T$ on
$\PP^s$ we will denote
\begin{eqnarray}
& & \mathcal P:=H^*_T\PP^s \\ & & \mathcal R=:\mathcal P\otimes \CC(\lambda)
\end{eqnarray}

Consider the diagonal action of $T=(\CC^*)^{s+1}$ on $\PP^s$ with weights $(-\lambda_0,...,-\lambda_s)$ i.e.
\begin{equation}
(t_0,t_1,...,t_s)\cdot (z_0,z_1,...,z_s)=(t_0^{-1}z_0,...,t_s^{-1}z_s).
\end{equation}
Then $\PP^s_T=\PP(\oplus_{i}\mathcal O(-\lambda_i))$. There is an
obvious lifting of the action of $T$ on the tautological line
bundle $\mathcal O(-1)$. It follows that $\mathcal O(k)$ is
equivariant for all $k$. Let $p=c_1^T(\mathcal O_{\PP^s}(1))$ be
the equivariant hyperplane class. We obtain $\mathcal
P=\CC[\lambda,p]/\prod_{i}(p-\lambda_i)$ and $\mathcal
R=\CC(\lambda)[p]/\prod_{i}(p-\lambda_i)$. The locus of the fixed
points consists of points $p_j$ for $j=0,1,...,s$ where $p_j$ is
the point whose $j$-th coordinate is $1$ and all the other ones
are $0$. On the level of the cohomology the map $i_j^*$ sends $p$
to $\lambda_j$. A basis for $\mathcal R$ as a
$\CC(\lambda)$-vector space is given by $\phi_j=\prod_{k\neq
j}(p-\lambda_k)$ for $j=0,1,...,s$. Also
$i_j^*(\phi_j)=\prod_{k\neq
j}(\lambda_j-\lambda_k)=\text{Euler}_T(N_j)$. The localization
theorem $2.2.1$ says that for any polynomial $F(p)\in
\CC(\lambda)[p]/\prod_{i=1}^{s}(p-\lambda_i)$
\begin{equation}
\int_{\PP^s_T}F(p)=\sum_{j}\frac{F(\lambda_j)}{\prod_{k\neq j}(\lambda_j-\lambda_k)}.
\end{equation}
By translating the target of a stable map we get an action of $T$
on ${\overline M}_{0,n}(\PP^s,d)$. In \cite{[24]} Kontsevich
identified the fixed point components of this action in terms of
decorated graphs. If $f:(C,x_1,...,x_n)\rightarrow \PP^s$ is a
fixed stable map then $f(C)$ is a fixed curve in $\PP^s$. The
marked points, collapsed components and nodes are mapped to the
fixed points $p_i$ of the $T$-action on $\PP^s$. A noncontracted
component must be mapped to a fixed line $\overline{p_ip_j}$ on
$\PP^s$. The only branch points are the two fixed points $p_i$ and
$p_j$ and the restriction of the map $f$ to this component is
determined by its degree. The graph $\Gamma$ corresponding to the
fixed point component containing such a map is constructed as
follows. The vertices correspond to the connected components of
$f^{-1}\{p_0,p_1,...,p_s\}$. The edges correspond to the
noncontracted components of the map. The graph is decorated as
follows. Edges are marked by the degree of the map on the
corresponding component, and vertices are marked by the fixed
point of $\PP^s$ where the corresponding component is mapped to.
To each vertex we associate a leg for each marked point that
belongs to the corresponding component. For a vertex $v$, let
$n(v)$ be the number of legs or edges incident to that vertex.
Also for an edge $e$ let $d_e$ be the degree of the stable map on
the corresponding component. Let
\begin{equation}
\overline {\mathcal M}_{\Gamma}:=\prod_{v}{\overline M}_{0,n(v)}.
\end{equation}
There is a finite group of automorphisms $\text{G}_{\Gamma}$
acting on ${\overline M}_{\Gamma}$ \cite{[9]}\cite{[16]}. The
order of the automorphism group $\text{G}_{\Gamma}$ is
\begin{equation}
a_{\Gamma}=\prod_{e}d_e\cdot |\text{Aut}(\Gamma)|. \label{eq: autorder}
\end{equation}
The fixed point component corresponding to the decorated graph $\Gamma$ is
\begin{equation}
\overline M_{\Gamma}={\overline {\mathcal M}}_{\Gamma}/{\text G}.
\end{equation}
Let $i_{\Gamma}:\overline  M_\Gamma\hookrightarrow \overline
{M}_{0,n}(\PP^s,d)$ be the inclusion of the  fixed point component
corresponding to $\Gamma$ and $N_{\Gamma}$ its normal bundle. This
bundle is $T$-equivariant. Let $\alpha$ be an equivariant
cohomology class in $H^*_T({\overline M}_{0,n}(\PP^s,d))$ and
$\alpha_{\Gamma}:=i^*_{\Gamma}(\alpha)$. Theorem $2.2.1$ says:
\begin{equation}
\int_{{\overline
M}_{0,n}(\PP^s,d)_T}\alpha=\sum_{\Gamma}\int_{({\overline
M}_{\Gamma})_T}\frac{\alpha_{\Gamma}}{a_{\Gamma}\text{Euler}_T(N_{\Gamma})}.
\end{equation}
Explicit formulas for $\text{Euler}_T(N_{\Gamma})$ in terms of
chern classes of cotangent line bundles in $H^*_T({\overline
M}_{\Gamma})$ have been found by Kontsevich in \cite{[24]}.

\subsection{\bf Linear and nonlinear sigma models for a projective space.}
Two compactifications of the space of
degree $d$ maps $\PP^1\rightarrow \PP^s$ will be very important in
this paper. $M_d:={\overline M}_{0,0}({\PP}^s\times{\PP}^1,(d,1))$
is called the degree $d$ {\it nonlinear sigma model} of $\PP^s$
and $N_d:=\PP(H^0(\PP^1,\mathcal O_{\PP^1}(d))^{s+1})$ is called
the degree $d$ {\it linear sigma model} of the projective space
$\PP^s$. An element in $H^0(\PP^1,\mathcal O_{\PP^1}(d))^{s+1}$ is
an $s+1$-tuple of degree $d$ homogeneous polynomials in two
variables $w_0$ and $w_1$. As a vector space, $H^0(\PP^1,\mathcal
O_{\PP^1}(d))^{s+1}$ is generated by the vectors
$v_{ir}=(0,...,0,w_0^rw_1^{d-r},0...,0)$ for $i=0,1,...,s$ and
$r=0,1,...,d$. The only nonzero component of $v_{ir}$ is the
$i$-th one.

The action of $T':=T\times {\bf{\CC}^*}$ in ${\PP}^s\times
{\PP}^1$ with weights $(-\lambda_0,...,-\lambda_s)$ in the $\PP^s$
factor and $(-\hbar,0)$ in the $\PP^1$ factor gives rise to an
action of $T'$ in $M_d$ by translation of maps. $T'$ also acts in
$N_d$ as follows: for $\bar{t}=(t_0,...,t_s)\in T$ and $t\in
{\bf{\CC}^*}$
\begin{equation}
(\bar{t},t)\cdot [P_0(w_0,w_1),...,P_s(w_0,w_1)]=[t_0P_0(tw_0,w_1),...,t_sP_s(tw_0,w_1)]
\end{equation}
There is a $T'$-equivariant morphism $\psi:M_d\mapsto N_d.$ Here
is a set-theoretical description of this map (for a proof that it
is a morphism see \cite{[15]} or \cite{[25]}). Let $q_i$ for
$i=1,2$ be the projection maps on $\PP^s\times \PP^1$. For a
stable map $(C,f)\in M_d$ let $C_0$ be the unique component of $C$
such that $q_2\circ f:C_0\rightarrow \PP^1$ is an isomorphism. Let
$C_1,...,C_n$ be the irreducible components of $C-C_0$ and $d_i$
the degree of the restriction of $q_1\circ f$ on $C_i$. Choose
coordinates on $C_0\cong{\PP^1}$ such that $q_2\circ
f(y_0,y_1)=(y_1,y_0)$. Let $\displaystyle{C_0\cap C_i=(a_i,b_i)}$
and $\displaystyle{q_1\circ f=[f_0:f_1:...:f_s]:C_0\mapsto
{\PP^s}}$. Then
\begin{equation}
\psi(C,f):=\prod_{i=1}^{n}(b_iw_0-a_iw_1)^{d_i}[f_0:f_1:...:f_s].
\end{equation}
Let $p_{ir}$ be the points of $N_d$ corresponding to the vectors
$v_{ir}$. The fixed point loci of the $T'$-action on $N_d$
consists of the points $p_{ir}$. We write $\kappa$ for the
equivariant hyperplane class of $N_d$. The restriction of $\kappa$
at the fixed point $p_{ir}$ is $\lambda_i+r\hbar$. The restriction
of the equivariant Euler class of the tangent space $TN_d$ at
$p_{ir}$ is \cite{[25]}
\begin{equation}
E_{ir}=\prod_{(j,t)\neq (i,r)}(\lambda_i-\lambda_j+r\hbar-t\hbar). \label{eq: EulerNd}
\end{equation}
Fixed point components of $M_d$ are obtained as follows. Let
$\Gamma^i_{d_j}$ be the graph of a $T$-fixed point component in
$\overline{M}_{0,1}(\PP^s,d_j)$ where the marking is mapped to
$p_i$ and $d_1+d_2=d$. Let $(d_1,d_2)$ be a partition of $d$. We
identify $\overline{M}_{\Gamma^i_{d_1}}\times
\overline{M}_{\Gamma^i_{d_2}}$ with a $T'$-fixed point component
$M^i_{d_1d_2}$ in $M_d$ in the following manner. Let
$(C_1,x_1,f_1)\in \overline{M}_{\Gamma^i_{d_1}}$ and
$(C_2,x_2,f_2)\in \overline{M}_{\Gamma^i_{d_2}}$. Let $C$ be the
nodal curve obtained by gluing $C_1$ with $\PP^1$ at $x_1$ and
$0\in \PP^1$ and $C_2$ with $\PP^1$ at $x_2$ and $\infty\in
\PP^1$. Let $f:C\rightarrow \PP^s\times \PP^1$ map $C_1$ to the
slice $\PP^s\times \infty$ by means of $f_1$ and $C_2$ to
$\PP^s\times 0$ by means of $f_2$. Finally $f$ maps $\PP^1$ to
$p_i\times \PP^1$ by permuting coordinates. $\psi$ maps
$M^i_{d_1d_2}$ to $p_{id_{2}}\in N_d$, hence the equivariant
restriction of $\psi^*(\kappa)$ in $M^i_{d_1d_2}$ is
$\lambda_i+d_2\hbar$. The normal bundle $N_{\Gamma^i_{d_1d_2}}$ of
this component in the above identification can be found by
splitting it in five pieces: smoothing the nodes $x_1$ and $x_2$
and  deforming the restriction of the map to $C_1,C_2,\PP^1$.
Using Kontsevich's calculations, Givental obtained \cite{[15]}
\[\frac{1}{\text{E}_T(N_{\Gamma^i_{d_1d_2}})}=\frac{1}{\prod_{k\neq i}(\lambda_i-\lambda_k)}\frac{1}{\text{E}_T(N_{\Gamma^i_{d_1}})}\frac{1}{\text{E}_T(N_{\Gamma^i_{d_2}})}\frac{e_1^*(\phi_i)}{-\hbar(-\hbar-c_1)}\frac{e_1^*(\phi_i)}{\hbar(\hbar-c_2)}\]
where $c_j , j=1,2$ is the first Chern class of the
cotangent line bundle on $\overline{M}_{\Gamma^i_{d_j}}$.

\section{A Gromov-Witten theory induced by a vector bundle}
\subsection{\bf The obstruction class of a concavex vector bundle.} The notion of concavex vector
bundle is due to Lian, Liu and Yau \cite{[25]} and is central to this work.
\begin{definition}
\begin{enumerate}
\item A line bundle $\mathcal L$ on $X$ is called {\bf convex} if $H^1(C,f^*(\mathcal L))=0$ for any genus zero stable map $(C,x_1,...x_n,f)$.
\item A line bundle $\mathcal L$ on $X$ is called {\bf concave} if $H^0(C,f^*(\mathcal L))=0$ for any nonconstant genus zero stable map $(C,x_1,...x_n,f)$.
\item A direct sum of convex and concave line bundles on $X$ is called a {\bf concavex} vector bundle.
\end{enumerate}
\end{definition}

A  concavex vector bundle $V$ in a projective space $\PP^s$ has
the form
\begin{equation}
V=V^+\oplus V^-=\left(\oplus_{i\in I}\mathcal O(k_i)\right)\oplus
\left(\oplus _{j\in J} \mathcal O(-l_j)\right)
\end{equation}
where $k_i$ and $l_j$ are positive numbers. Denote
$\EE^+:=\text{E}(V^+)$ and $\EE^-:=\text{E}(V^-)$.

Let $d>0$. Consider the following diagram
\[ \begin{CD}
{{\overline M}_{0,n+1} ({\PP}^s,d)}@>e_{n+1}>>{\PP}^s \\
@VV \pi_{n+1} V \\
{{\overline M}_{0,n} ({\PP}^s,d)}
\end{CD} \]
Since $V$ is concavex, the sheaf
\begin{equation}
V_d:=V_d^+\oplus V_d^-= {\pi_{n+1}}_*e_{n+1}^*(V^+)\oplus
R^1{\pi_{n+1}}_*e_{n+1}^*(V^-)
\end{equation}
is locally free.
\begin{definition} The obstruction class corresponding to $V$ is defined to be:
\begin{equation}
\EE_d:=\text
{E}(V_d)=\text{E}(V^+_d)\text{E}(V^-_d):=\EE_d^+\EE_d^-
\end{equation}
\end{definition}
For a $T$-action on $\PP^s$ that lifts to a linear action on the
fibers of $V=V^+\oplus V^-$, let $E^+:=\text{E}_T(V^+)$ and
$E^-:=\text{E}_T(V^-)$. Assume that $E^-$ is invertible.
\begin{definition}
The modified equivariant integral  $\omega_V: \mathcal
R\rightarrow \CC(\lambda)$ corresponding to $V$ is defined as
follows
\begin{equation}
\omega_V(\alpha):=\int_{\PP_T^m}\alpha\cup \frac{E^+}{E^-}.
\end{equation}
\end{definition}
Consider the trivial action of $T=(\oplus_{i\in I}C^*)\oplus
(\oplus_{j\in J}C^*)$ on $\PP^s$. In this case
$\PP^s_T=\PP^s\times (\oplus_{i\in I}\PP^{\infty})\times
(\oplus_{j\in J}\PP^{\infty})$ and ${{\overline M}_{0,n}
({\PP}^s,d)}_{T}={\overline M}_{0,n} ({\PP}^s,d) \times
(\oplus_{i\in I}\PP^{\infty})\times (\oplus_{j\in
J}\PP^{\infty})$. It follows that $\mathcal
P=H^*(\PP^s,\CC[\lambda])$ and $\mathcal
R=H^*(\PP^s,\CC(\lambda))$. Let $p$ denotes the equivariant
hyperplane class. The $T$-action lifts to a linear action on the
fibers of $V$ with weights $((-\lambda_i)_{i\in
I},(-\lambda_j)_{j\in J})$. Let $q_i$ and $q_j$ denote the
projection maps on ${{\overline M}_{0,n} ({\PP}^s,d)}_{T}$. Both
$V^+_d$ and $V^-_d$ are $T$-equivariant bundles and
\begin{eqnarray}
& & (V^+_d)_T=V^+_d\otimes (\oplus_{i\in I}q_i^*{\mathcal
O}_{\PP^\infty}(-\lambda_i)) \nonumber \\ & &
(V^-_d)_T=V^-_d\otimes (\oplus_{j\in J}q_j^*{\mathcal
O}_{\PP^\infty}(-\lambda_j))
\end{eqnarray}
The equivariant obstruction class is
\begin{equation}
E_d:=\text{E}_T(V_d)=\text{E}_T(V^+_d)\text{E}_T(V^-_d)=E^+_dE^-_d.
\label{eq: euler}
\end{equation}
The modified equivariant integral for the trivial action of $T$ on
$\PP^s$ gives rise to a modified perfect pairing in $\mathcal R$
\[{\langle a,b\rangle}_V:=\omega_V(a\cup b).\]
Let $T_0=1,T_1=p,...,T^s=p^s$ be a basis of $\mathcal R$ as a
$\CC(\lambda)$-vector space. The intersection matrix
$(g_{rt}):=({\langle T_r,T_t\rangle}_V)$ has an inverse
$(g^{rt})$. Let $T^i=\sum_{j=0}^{s}g^{ij}T_j$ be the dual basis
with respect to this pairing. Clearly
\begin{equation}
T^i=T_{m-i}\cdot \left(\frac{E^-}{E^+}\right). \label{eq: basis}
\end{equation}
This implies that in $H^*(\PP^s\times \PP^s)\otimes \CC(\lambda)$
we have
\begin{equation}
\sum_{i=1}^{s}T_i\otimes T^i=\Delta \cdot \left(1\otimes
\frac{E^-}{E^+}\right)
\end{equation}
where $\Delta=\sum_{i=0}^{s}T_i\otimes T_{s-i}$ is the class of
the diagonal in $\PP^s\times \PP^s$.

Recall that the morphism $\pi_k:{\overline M}_{0,n}
({\PP}^s,d)\rightarrow {\overline M}_{0,n-1} ({\PP}^s,d)$ forgets
the $k$-th marked point.
\begin{lemma} \label{lemma: pullback}
${\pi_k}^*(E_d)=E_d$ and ${\pi_k}^*(\EE_d)=\EE_d.$
\end{lemma}
{\bf Proof.} For simplicity we will consider the case $V=\mathcal
O(k)\oplus \mathcal O(-l)$ and $k=n$. The general case is similar.
\ Let $M_k=\overline M_{0,k} ({\PP}^s,d)$ and
$M_{n,n}=M_n\times_{M_{n-1}}M_n$. Consider the following
equivariant commutative diagram:
\[
\xymatrix
{& M_{n+1} \ar[ddl]_{\pi_{n+1}} \ar[d]_{\mu} \ar[drr]^{e_{n+1}} \ar[ddr]^{\pi_n} & & \\
& M_{n,n} \ar[dl]^{\beta} \ar[dr]
_{\alpha} & & \PP^s \\
M_n \ar[dr]_{\pi_n} & & M_n \ar[ur]_{e_n} \ar[dl]^{\pi_n} &  \\ &
M_{n-1} & & }
\]
We compute:
\begin{equation}
{\pi_{n+1}}_*e_{n+1}^*{\mathcal
O}(k)={\pi_{n+1}}_*\pi_n^*{e_n}^*{\mathcal
O}(k)=\beta_*\mu_*\mu^*\alpha^*{e_n}^*{\mathcal O}(k). \label{eq:
oneplus}
\end{equation}
By the projection formula
\begin{equation}
\mu_*\mu^*\alpha^*{e_n}^*{\mathcal O}(k)=\alpha^*{e_n}^*{\mathcal
O}(k)\otimes \mu_*(\mathcal O_{M_{n+1}}).
\end{equation}
Since the map $\mu$ is birational and $M_{n+1}$ is normal $\mu_*(\mathcal O_{M_{n+1}})=\mathcal O_{M_{n,n}}$
hence
\begin{equation}
\mu_*\mu^*\alpha^*{e_n}^*{\mathcal O}(k)=\alpha^*{e_n}^*{\mathcal
O}(k). \label{eq: projection}
\end{equation}
Substituting in (\ref{eq: oneplus}) and applying base extension
properties ($\pi_n$ is flat) yields
\begin{equation}
{\pi_{n+1}}_*e_{n+1}^*{\mathcal
O}(k)=\beta_*\alpha^*{e_n}^*{\mathcal
O}(k)={\pi_n}^*({\pi_n}_*{e_n}^*{\mathcal O}(k)).
\end{equation}
For the case of a negative line bundle we have
\begin{equation}
R^1{\pi_{n+1}}_*e_{n+1}^*{\mathcal
O}(-l)=R^1{\pi_{n+1}}_*{\pi_n}^*{e_n}^*{\mathcal
O}(-l)=R^1{\pi_{n+1}}_*\mu^*\alpha^*{e_n}^*{\mathcal O}(-l).
\end{equation}
We now use the spectral sequence
\begin{equation}
R^p{\beta}_*(R^q\mu_*\mathcal F)\Longrightarrow
R^{p+q}{\pi_{n+1}}_*{\mathcal F}
\end{equation}
where $\mathcal F$ is a sheaf of ${\mathcal O}_{\mathcal M_{n+1}}$-modules. The map $\mu$ is birational. If we think
of $M_n$ as the universal map of $M_{n-1}$, then the map $\mu$ has nontrivial fibers only over pairs of stable maps
in $M_n$ that represent the same special point (i.e. node or marked point) of a stable map in $M_{n-1}$. These
nontrivial fibers are isomorphic to $\PP^1$. Since $\mathcal F=e_{n+1}^*\mathcal O(-l)$ we obtain $R^q\mu_*\mathcal F=0$
for $q>0$. It follows that this spectral sequence degenerates, giving
\begin{equation}
R^1{\pi_{n+1}}_*e_{n+1}^*{\mathcal
O}(-l)=R^1\beta_*\mu_*\mu^*\alpha^*{e_n}^*{\mathcal O}(-l).
\end{equation}
Now we proceed as in (\ref{eq: projection}) to conclude
\begin{equation}
R^1{\pi_{n+1}}_*e_{n+1}^*{\mathcal
O}(-l)={\pi_n}^*(R^1{\pi_n}_*{e_n}^*{\mathcal O}(-l)).
\end{equation}
The lemma is proven.$\dagger$
\begin{remark}
The previous lemma justifies the omission of $n$ from the notation
of the obstruction class.
\end{remark}

\subsection{\bf Modified equivariant correlators and quantum cohomology.} Let $\gamma_i\in \mathcal R$ for $i=1,...,n$
and $d>0$. Introduce the following modified equivariant Gromow-Witten invariants:
\begin{equation}
\tilde{I_d}(\gamma_1,...,\gamma_n):=\int_{{\overline M_{0,n}
({\PP}^m,d)}_T} e_1^*(\gamma_1)\cup...\cup e_n^*(\gamma_n)\cup
E_d\in \CC(\lambda) \label{def: GWM}
\end{equation}
Now $\overline{M}_{0,n} ({\PP}^s,0)={\overline M}_{0,n}\times
\PP^s$ and all the evaluation maps equal the projection $q_2$ to
the second factor. The integrals in this case are defined as
follows
\begin{equation}
\tilde{I_0}(\gamma_1,...,\gamma_n):=\int_{\overline{M}_{0,n}
({\PP}^s,0)}e_1^*(\gamma_1)\cup...\cup e_n^*(\gamma_n)\cup
q_2^*(\text{E}(V))\in \CC(\lambda)
\end{equation}
The modified equivariant gravitational descendants are defined
similarly to Gromov-Witten invariants:
\begin{equation}
\tilde{I}_d(\tau_{k_1}\gamma_1,...,\tau_{k_n}\gamma_n):=
\int_{{{\overline M}_{0,n} ({\PP}^s,d)}_T}c_1^{k_1}({\mathcal
L}_1)\cup e_1^*(\gamma_1)\cup...\cup c_1^{k_n}(\mathcal L_n)\cup
e_n^*(\gamma_n) \cup E_d.
\end{equation}

Lemma (\ref{lemma: pullback}) is essential in proving that the
modified correlators satisfy the same properties, such as
fundamental class property, divisor property, point mapping axiom
etc., that the usual Gromov-Witten invariants do. The proofs are
similar to the ones in pure Gromov-Witten theory. As an
illustration, we prove one of these properties.

{\bf Fundamental class property}. Let $\gamma_n=1$ and $d\neq 0$.
The forgetful morphism $\pi_n:{\overline M}_{0,n}
({\PP}^s,d)\rightarrow {\overline M}_{0,n-1} ({\PP}^s,d)$ is
equivariant. Using Lemma~\ref{lemma: pullback} we obtain:
\[e_1^*(\gamma_1)\cup...\cup e_{n-1}^*(\gamma_{n-1})\cup e_n^*(1)\cup E_d=\pi^*(e_1^*(\gamma_1)\cup...\cup e_{n-1}^*(\gamma_{n-1})\cup E_d).\] Therefore:
\[\tilde{I_d}(\gamma_1,...,\gamma_{n-1},1)=\int_{{\overline M}_{0,n} ({\PP}^s,d)}\pi^*(e_1^*(\gamma_1)\cup...\cup e_{n-1}^*(\gamma_{n-1})\cup E_d)=\]
\[=\int_{{\pi_n}_*({\overline M}_{0,n} ({\PP}^s,d))}e_1^*(\gamma_1)\cup...\cup e_{n-1}^*(\gamma_{n-1})\cup E_d=0.\]
The last equality is because the fibers of $\pi_n$ are positive
dimensional. If $d=0$, by the point mapping property we know that
the integral is zero unless $n=3$. In that case:
$\tilde{I_0}(\gamma_1,\gamma_2,1)=\langle
\gamma_1,\gamma_2\rangle$.$\dagger$

We will now prove a technical lemma which will be very useful
later. Let $A\cup B$ be a partition of the set of markings and
$d=d_1+d_2$. Let $D=D(A,B,d_1,d_2)$ be the closure in $\overline
M_{0,n}(\PP^s,d)$ of stable maps of the following type. The source
curve is a union $C=C_1\cup C_2$ of two lines meeting at a node
$x$. The marked points corresponding to $A$ are on $C_1$ and those
corresponding to $B$ are on $C_2$. The restriction of the map $f$
on $C_i$ has degree $d_i$ for $i=1,2$. $D$ is a boundary divisor
in ${\overline M}_{0,n}(\PP^s,d)$. Let $M_1:=\overline
M_{0,|A|+1}(\PP^s,d_1)$ and $M_2:=\overline
M_{0,|B|+1}(\PP^s,d_2)$. Let $e_x$ and $\tilde{e_x}$ be the
evaluation maps at the additional marking in $M_1$ and $M_2$ and
$\mu:=(e_x,\tilde{e_x})$. The boundary divisor $D$ is obtained
from the following fibre diagram
\[
\begin{CD}
D                      @>\iota>>                M_1\times  M_2 \\
@VV\nu V                                         @VV\mu V\\
\PP^s                 @>\delta>>              \PP^s\times \PP^s\\
\end{CD}
\]
where $\nu$ is the ``evaluation map at the node $x$'' and $\delta$
is the diagonal map.
\begin{lemma} \label{lemma: split}
For any classes
$\gamma_1,...,\gamma_n$ in $\mathcal R$:
\[\int_{D}\prod_{i=1}^{n}e_i^*(\gamma_i)E_d=
\sum_{a=0}^{s} \left(\int_{M_1}\prod_{i\in
A}e_i^*(\gamma_i)e_x^*(T_a)E_{d_1}\right)\times
\left(\int_{M_2}\prod_{j\in
B}e_j^*(\gamma_j)\tilde{e}_x^*(T^a)E_{d_2}\right).\]
\end{lemma}
{\bf Proof.} This lemma is the analogue of the Lemma $16$ in
\cite{[13]}. The proof needs a minor modification. Let
$\alpha:D\rightarrow {\overline M}_{0,n} (\PP^s,d)$. Consider the
normalization sequence at $x$:
\begin{equation}
0\rightarrow {\mathcal O}_C\rightarrow {\mathcal O}_{C'}\oplus {\mathcal O}_{C''}\rightarrow {\mathcal O}_x\rightarrow 0.   \label{seq:norm}
\end{equation}
Twisting it by $f^*(V^+)$ and $f^*(V^-)$ and taking the cohomology sequence yields the following identities on $D$:
\begin{equation}
\alpha^*(E^+_d)\nu^*(E^+)=\iota^*(E^+_{d_1}\times E^+_{d_2}).
\label{eq: -E}
\end{equation}
and
\begin{equation}
\alpha^*(E^-_d)=\iota^*(E^-_{d_1}\times E^-_{d_2})\nu^*(E^-).
\label{eq: +E}
\end{equation}
By combining equations (\ref{eq: -E}) and (\ref{eq: +E}) we obtain the restriction of $E_d$ in the divisor $D$:
\begin{equation}
\alpha^*(E_d)=\iota^*(E_{d_1}\times
E_{d_2})\nu^*\left(\frac{E^-}{E^+}\right). \label{eq: EinD}
\end{equation}
Using formula $(23)$ we obtain
\begin{equation}
\iota_*\nu^*\left(\frac{E^-}{E^+}\right)=\mu^*\left(1\otimes
\frac{E^-}{E^+}\right)\mu^*(\Delta).
\end{equation}
Therefore
\begin{eqnarray}
& & \int_{D}\prod_{i=1}^{n}e_i^*(\gamma_i)\cup E_d=\int_{M_1\times
M_2}\prod_{i=1}^{n}e_i^*(\gamma_i)\cup E_{d_1}\cup E_{d_2}\cup
\mu^*\left(1\otimes \frac{E^-}{E^+}\right)\cup
\nonumber\\
& & \cup \mu^*(\Delta)= \int_{M_1\times
M_2}\prod_{i=1}^{n}e_i^*(\gamma_i)\cup E_{d_1}\cup E_{d_2}\cup
\mu^*\left(\sum_{a}T_a\otimes T^a\right)=
\nonumber \\
& & \sum_{a=0}^{m}
\left(\int_{M_1}\prod_{i=1}^{n_1}e_i^*(\gamma_i)\cup
e_x^*(T_a)\cup E_{d_1}\right)\times
\left(\int_{M_2}\prod_{j=1}^{n_2}e_j^*(\gamma_j)\cup
\tilde{e}_x^*(T^a)\cup E_{d_2}\right).
\end{eqnarray}
The lemma is proven.$\dagger$

The same proof can be used to show that the previous splitting
lemma is true for gravitational descendants as well.
\begin{corollary}
The following modified topological recursion relations hold:
\begin{eqnarray}
& & \tilde{I}_d(\tau_{k_1+1}\gamma_1,\tau_{k_2}\gamma_2,\tau_{k_3}\gamma_3,
\prod_{i=4}^{n}\tau_{s_i}\omega_i)= \nonumber \\ & & \sum_{}\tilde{I}_{d_1}(\tau_{k_1}\gamma_1,\prod_{i\in I_1}\tau_{s_i}\omega_i,T_a)\tilde{I}_{d_2}(T^a,\tau_{k_2}\gamma_2,\tau_{k_3}\gamma_3,\prod_{i\in I_2}\tau_{s_i}\omega_i) \label{eq: toprec}
\end{eqnarray} where the sum is over all splittings $d_1+d_2=d$ and partitions $I_1\cup I_2=\{4,...,n\}$ and over all indices $a$.
\end{corollary}
{\bf Proof.} Let $A$ and $B$ be two disjoint subsets of
$\{1,2,...,n\}$. We will denote by $D(A,B)$ the sum of boundary
divisors $D(E,F,d_1,d_2)$ such that $E, F$ is a partition of
$\{1,2,...,n\}$ and $A\subset E$ and $B\subset F$ and $d_1+d_2=d$.
The notation $D(A,B)$ does not reflect neither the number $n$ of
marked points nor the degree $d$ of the maps but they will be
clear from the context. Consider the morphism
$\pi:\overline{M}_{0,n}(\PP^s,d)\rightarrow \overline M_{0,3}$
that forgets the map and all but the first $3$ markings. Since
$\overline M_{0,3}$ is a point, the cotangent line bundle at the
first marking is trivial. But $\pi^*(\mathcal L_1)=\mathcal
L_1-D(\{1\},\{2,3\})$ therefore $\mathcal L_1=D(\{1\},\{2,3\})$ in
$\overline{M}_{0,n}(\PP^s,d)$. Multiply both sides of the previous
equation by $\prod_{j=1}^{3}c_1(\mathcal L_j)^{k_j}\cup
e_j^*(\gamma_j)\cup \prod_{i=4}^{n}c_1(\mathcal L_i)^{s_i}\cup
e_i^*(\omega_i)\cup E_d$ and integrate. The corollary follows from
the splitting lemma for gravitational descendents.$\dagger$

In the process of finding solutions to the WDVV equations,
Kontsevich suggested the following modified equivariant
Gromov-Witten potential
\begin{equation}
\tilde{\Phi}(t_0,t_1,...,t_m):=\sum_{n\geq 3}\sum_{d\geq 0}\frac{1}{n!}\tilde{I_d}(\gamma^{\otimes n}) \label{def: QPM}
\end{equation} where $\gamma=t_0+t_1p+...+t_sp^s$ and $t_i\in
\CC(\lambda)$. Let
$\displaystyle{\tilde{\Phi}_{ijk}=\frac{\partial^3\tilde{\Phi}}{\partial
t_i \partial t_j \partial t_k}}$.
\begin{definition} The modified, equivariant quantum product on $\mathcal R$ is defined to be the linear extension of
\begin{equation}
T_i*_V T_j:=\sum_{k=0}^{m}\tilde{\Phi}_{ijk}T^k.
\end{equation}
\end{definition}
\begin{theorem} \label{theorem: ring}
$QH^*_V\PP^s_T:=({\mathcal R},*_V)$ is a commutative, associative
algebra with unit $T_0$.
\end{theorem}
{\bf Proof.} A simple calculation shows that:
\[\tilde{\Phi}_{ijk}=\sum_{n\geq 0}\sum_{d\geq
0}\frac{1}{n!}\tilde{I_d}(T_i,T_j,T_k,\gamma^{\otimes n}).\] The
commutativity of the modified, equivariant quantum product follows
from the symmetry of the new integrals. $T_0$ is the unit due to
the fundamental class property for the modified Gromov-Witten
invariants. To proving the associativity we proceed as in Theorem
4 in \cite {[12]}. Let
$\displaystyle{\tilde{\Phi}_{ijk}=\frac{\partial^3\tilde{\Phi}}{\partial
t_i \partial t_j \partial t_k}}.$ We compute
\[(T_i*_VT_j)*_VT_k=\sum_{}\sum{}\tilde{\Phi}_{ije}g^{ef}\tilde{\Phi}_{fkl}g^{ld}T_d\]
\[T_i*_V(T_j*_VT_k)=\sum_{}\sum_{}\tilde{\Phi}_{jke}g^{ef}\tilde{\Phi}_{fil}g^{ld}T_d.\]
Since the matrix $(g^{ld})$ is nonsingular, $(T_i*_VT_j)*_VT_k=T_i*_V(T_j*_VT_k)$ is equivalent to
\begin{equation}
\sum_{e,f}\tilde{\Phi}_{ije}g^{ef}\tilde{\Phi}_{fkl}=\sum_{e,f}\tilde{\Phi}_{jke}g^{ef}\tilde{\Phi}_{fil}. \label{eq: WDVV}
\end{equation}
Equation (\ref{eq: WDVV}) is the WDVV equation for the modified potential $\tilde{\Phi}$. To prove this equation
let $q,r,s,t$ be four different integers in $\{1,2,...n\}$. There exists an equivariant morphism:
\[\pi:{\overline M}_{0,n}(\PP^s,d)\rightarrow {\overline M}_{0,4}=\PP^1\] that forgets the map and all the marked
points but $q,r,s,t$.  Obviously the divisors $D(\{q,r\},\{s,t\})$
and $D(\{q,s\},\{r,t\})$ are linearly equivalent in $\overline
M_{0,4}$ hence, via the pullback $\pi^*$, they are linearly
equivalent in ${\overline M}_{0,n}(\PP^s,d)$. Now integrate the
class
\[\prod_{i=1}^{n-4}(e_i^*(\gamma))\cup e_{n-3}^*(T_i)\cup e_{n-2}^*(T_j)\cup e_{n-1}^*(T_k)\cup e_n^*(T_l)\cup E_d\]
over $D(\{q,r\},\{s,t\})$ and use Lemma~\ref{lemma: split} to
obtain WDVV equation hence the associativity.$\dagger$

If we restrict $\tilde{\Phi}_{ijk}$ to the divisor classes
$\gamma=tp$, and use the divisor property for the modified
Gromov-Witten invariants, we obtain the {\it small product}:
\begin{equation}
T_i*_V T_j:=T_i\cup T_j+\sum_{d>0}q^d\sum_{k=0}^{m}\tilde{I_d}(T_i,T_j,T_k)T^k.  \label{eq: smq}
\end{equation}
Here $q=e^t$. We extend this product to $\mathcal
R\otimes_{\CC}\CC[[q]]$ to obtain the small equivariant quantum
cohomology ring $SQH^*_V{\PP^s}_T$. We will use $*_V$ to denote
both the small and the big quantum product. The difference will be
clear from the context.
\begin{remark} \label{eq: Remark}
\begin{itemize}
\item Equation (\ref{eq: EinD}) and Lemma~\ref{lemma: pullback} are the basis for building a modified
equivariant Gromow-Witten theory similar to  pure Gromov-Witten theory.
\item One can see from (\ref{eq: basis}) that the only potential problem with the existence of the nonequivariant
limit of (\ref{eq: smq}) is the presence of $E^+$ in the denominator of $T^k$.
Hence if $V=V^-$ is a pure negative line bundle the nonequivariant limit of this product exists.
An example of this situation is treated in the last section.
\end{itemize}
\end{remark}
\section{A $\mathcal D$-module structure induced by $V$}
\subsection{Equivariant quantum differential equations.}
Recall from section $2.3$ the generator $\hbar$ of $H^2(BC^*)$.
Consider the system of first order differential equations on the
modified, big quantum cohomology ring $QH^*_V(\PP^s_T)$
\begin{equation}
\hbar\frac{\partial}{\partial t_i}=T_i*_V : i=1,...,m.
\end{equation}
\begin{theorem}
The space of solutions of these equations has the following basis:
\begin{eqnarray}
& &
s_a=T_a+\sum_{j=0}^{s}\sum_{n=0}^{\infty}\sum_{d=0}^{\infty}\sum_{k=0}^{\infty}\frac{\hbar^{-(k+1)}}{n!}\tilde{I}_d(\tau_kT_a,T_j,\gamma^{\otimes
n})T^j= \nonumber \\ & &
T_a+\sum_{j=0}^{m}\sum_{d=0}^{\infty}\sum_{n=0}^{\infty}\frac{1}{n!}\tilde{I}_d\left(\frac{T_a}{\hbar-c},T_j,\gamma^{\otimes
n}\right)T^j
\end{eqnarray}
where $c$ is a formal symbol that stands for ${c_1}^T(\mathcal
L_1)$ and $\displaystyle{\frac{T_a}{\hbar-c}}$ should be expanded
in powers of $\displaystyle{\frac{c}{\hbar}}$.
\end{theorem}
{\bf Proof.} On one hand
\begin{equation}
\hbar\frac{\partial s_a}{\partial t_i}=\sum_{j=0}^{m}\sum_{n=0}^{\infty}\sum_{d=0}^{\infty}\sum_{k=0}^{\infty}\frac{\hbar^{-k}}{n!}\tilde{I}_d(\tau_kT_a,T_j,T_i,\gamma^{\otimes n})T^j.
\end{equation}
On the other hand
\begin{eqnarray}
& & T_i*s_a=T_i*T_a+\sum_{j=0}^{m}\sum_{n=0}^{\infty}\sum_{d=0}^{\infty}\sum_{k=0}^{\infty}\frac{\hbar^{-(k+1)}}{{n_1}!}\tilde{I}_d(\tau_kT_a,T_j,\gamma^{\otimes n})(T_i*T^j)= \nonumber \\ & & \sum_{n,d,e}\frac{1}{n!}\tilde{I}_d(T_i,T_a,T_e,\gamma^{\otimes n})T^e+\sum_{j=0}^{m}\sum_{n=0}^{\infty}\sum_{k=0}^{\infty}\sum_{d_1}\frac{\hbar^{-(k+1)}}{{n}!}\tilde{I}_{d_1}(\tau_kT_a,T_j,\gamma^{\otimes n}) \nonumber \\ & & \times \sum_{m,d_2,e}\frac{1}{{m}!}\tilde{I}_{d_2}(T_i,T^j,T_e,\gamma^{\otimes m})T^e.
\end{eqnarray}
The theorem follows from the topological recursion relations
(\ref{eq: toprec}).$\dagger$

Restrictions $\tilde{s}_a$ of the sections $s_a$ to $\gamma\in
H^0(\PP^m)\oplus H^2(\PP^m)$ are solutions of
\[\hbar\frac{\partial}{\partial t_i}=T_i*_V : i=0,1.\]  Repeated
use of the divisor axiom yields
\begin{equation}
\tilde{s}_a=e^{\frac{t_0+pt_1}{\hbar}}\cup
T_a+\sum_{d=1}^{\infty}\sum_{j=0}^{m}q^d\tilde{I}_d\left(\frac{e^{\frac{t_0+pt_1}{\hbar}}\cup
T_a}{\hbar-c},T_j\right)T^j,  \label{eq: sec}
\end{equation}
where $q:=e^{t_1}$
\begin{definition}
The module of differential operators that annihilate
${\langle\tilde{s}_a,1\rangle}_V$ for all $a$ is called the
modified equivariant $\mathcal D$-module of $\PP^s$ induced by
$V$.
\end{definition}
This module is generated by the following $\mathcal
R[[t_0,t_1,q]]$-valued function
\begin{equation}
\tilde{J}_V=\sum_{a=0}^{s}{\langle\tilde{s}_a,1\rangle}_VT^a.
\label{eq: Jdef}
\end{equation}
Recall that $e_1: \overline M_{0,2}(\PP^s,d)\rightarrow \PP^s$ is
the evaluation map at the first marked point and $c$ is the chern
class of the cotangent line bundle at the first marked point.
Substituting (\ref{eq: sec}) in (\ref{eq: Jdef}) and using the
projection formula we obtain:
\begin{eqnarray}
& & \tilde{J}_V=\text{exp}\left(\frac{t_0+pt_1}{\hbar}\right)\cdot \\
& & \cdot
\left(1+\sum_{d>0}q^dPD^{-1}\left({e_1}_*\left(\left(\frac{E_d}{\hbar-c}\cap
[\overline M_{0,2}(\PP^s,d)] \right)\right)\cup
\left(\frac{E^-}{E^+}\right)\right)\right) \nonumber.
\end{eqnarray}
In the above expression $PD:H^*(\overline
M_{0,2}(\PP^s,d))\rightarrow H_{s+d+sd-1-*}\overline
M_{0,2}(\PP^s,d))$ is the Poincar\'e duality isomorphism.

It will be convenient for us to work with the moduli space of one
pointed stable maps. To that end we note that
\begin{equation}
{e_1}_*\left(\frac{E_d}{\hbar-c}\cap [\overline M_{0,2}(\PP^s,d)]
\right)={e_1}_*\left(\frac{E_d}{\hbar(\hbar-c)}\cap [\overline
M_{0,1}(\PP^s,d)] \right).
\end{equation}
This identity follows easily from the fact that if $\pi_2:
\overline M_{0,2}(\PP^s,d)\rightarrow \overline M_{0,1}(\PP^s,d)$
forgets the second marked point and $D$ is the image of the
universal section of $\pi$ induced by the marked point, then
$c={\pi_2}^*(c)+D$ and $\EE_d={\pi_2}^*(\EE_d)$.

The final expression for $\tilde{J}_V$ is
\begin{eqnarray}
& & \tilde{J}_V=\text{exp}\left(\frac{t_0+pt_1}{\hbar}\right) \cdot \\
& & \cdot
\left(1+\sum_{d>0}q^dPD^{-1}\left({e_1}_*\left(\frac{E_d}{\hbar(\hbar-c)}\cap
[\overline M_{0,1}(\PP^s,d)] \right)\right)\cup
\left(\frac{E^-}{E^+}\right)\right). \nonumber
\end{eqnarray}
From this presentation we see that the presence of the equivariant
class $E^+$ in the denominator of $\tilde{J}_V$ is a potential
problem for the existence of the nonequivariant limit.
\begin{lemma}
$\tilde{J}_V\in \mathcal P[[q]]$ therefore it has a nonequivariant limit.
\end{lemma}
{\bf Proof.} Let $V'_d$ be the subbundle of $V^+_d$ whose fiber
consists of those sections of $H^0(C,f^*(V^+))$ that vanish at the
marked point. Let $E'_d:=c^T_{\text{top}}(V'_d)$. There is an
exact sequence of equivariant bundles on ${\overline
M}_{0,1}(\PP^s,d)$:
\[0\rightarrow  V'_d\rightarrow V^+_d\rightarrow {e_1}^*(V^+)\rightarrow 0.\]

Taking the top chern classes we obtain
\begin{equation}
E^+_d=E'_d\cdot e_1^*(E^+). \label{eq: Eulersplit}
\end{equation}

We compute
\begin{eqnarray}
 & & PD^{-1}{e_1}_*\left(\frac{E_d}{\hbar(\hbar-c)}\cap
[\overline M_{0,1}(\PP^s,d)]\right)= \nonumber \\
& &
PD^{-1}\left({e_1}_*\left(\frac{E'_de_1^*(E^+)E_d^-}{\hbar(\hbar-c)}\cap
[\overline
M_{0,1}(\PP^s,d)] \right)\right) = \nonumber \\
& &  PD^{-1}\left(E^+\cap
{e_1}_*\left(\frac{E'_dE_d^-}{\hbar(\hbar-c)}\cap [\overline
M_{0,1}(\PP^s,d)]
\right)\right)= \nonumber \\
& & E^+\cup
PD^{-1}\left({e_1}_*\left(\frac{E'_dE_d^-}{\hbar(\hbar-c)}\cap
[\overline M_{0,1}(\PP^s,d)] \right)\right).
\end{eqnarray}
Therefore
\begin{eqnarray}
& & \tilde
{J}_V=\text{exp}\left(\frac{t_0+pt_1}{\hbar}\right)\cdot \nonumber \\
& & \cdot
\left(1+\sum_{d>0}q^dPD^{-1}{e_1}_*\left(\frac{E'_dE_d^-}{\hbar(\hbar-c)}\cap
[\overline M_{0,1}(\PP^s,d)] \right)\cup E^-\right).
\end{eqnarray}
It is now visible from this presentation that $\tilde{J}_V\in
\mathcal P[[q]]$ and
\begin{eqnarray}
& &J_V:=\lim_{\lambda\rightarrow 0}\tilde{J}_V= \nonumber \\
& &
\text{exp}\left(\frac{t_0+Ht_1}{\hbar}\right)\left(1+\sum_{d>0}q^dPD^{-1}{e_1}_*\left(\frac{\EE'_d\EE_d^-
}{\hbar(\hbar-c)}\cap [\overline M_{0,1}(\PP^s,d)]\right)\cup
\EE^-\right).
\end{eqnarray}
The lemma is proven.$\dagger$
\subsection{\bf A local property of the $J$-function} Let $Y$ be a smooth projective variety and
$j:\PP^s\hookrightarrow Y$ an embedding. Suppose that $\mathcal
N_{\PP^s/Y}=V^-={\mathcal O}(-l)$ for some $l>0$. Let $C$ be a
curve in $\PP^s$. The map $j$ gives rise to an embedding
\begin{eqnarray}
\overline M_{0,n}(\PP^s,[C])\hookrightarrow \overline
M_{0,n}(Y,j_*([C]))
\end{eqnarray}
\begin{lemma} \label{lemma: eqmod}
Let $C$ be a degree $d$ rational curve in $\PP^s$. Then $\overline
M_{0,n}(\PP^s,d)=\overline M_{0,n}(Y,j_*([C])).$
\end{lemma}
{\bf Proof.} Let $(C',x_1,...,x_n,f)\in {\overline
M}_{0,n}(Y,j_*([C]))$ and $f(C')=C_1\cup C_2\cup...\cup C_p$ be
the irreducible decomposition. Then
$d[\text{line}]=[C_1]+...[C_p]$. Let $I_1=\{i:
C_i\subset{\PP}^s\}$ and $I_2=\{1,2,...,n\}-I_1$. Assume that
$I_2$ is nonempty. If
\[d[\text{line}]-\sum_{i\in I_1}[C_i]\] has nonpositive degree in
$\PP^s$, we intersect with an ample divisor in $Y$ to see that
\[d[\text{line}]-\sum_{i\in I_1}[C_i]=\sum_{i\in I_2}[C_i]\]
is impossible. Otherwise, we intersect with $[\PP^s]$ to get the
same contradiction. Hence $I_2$ is empty and all the curves $C_i$
lie in $\PP^s$. It follows that $f$ factors through $\PP^s$ and
therefore $(C',x_1,...x_n,f)\in \overline M_{0,n}(\PP^s,d)$. On
the other hand $\overline M_{0,n}(\PP^s,d)$ is a component of
$\overline M_{0,n}(Y,j_*([C]))$ (see for example section $7.4.4$
in \cite{[9]}). These two arguments imply the lemma.$\dagger$

Denote ${\overline M}_{0,n}(Y,d):=\overline M_{0,n}(Y,j_*([C]))$,
where $C$ is any rational curve of degree $d$ in $\PP^s$. The
following Lemma is a special case of a conjecture by Cox, Katz and
Lee in \cite{[10]} which was proved in \cite{[23]}.
\begin{lemma} \label{lemma: refeuler}
$[\overline M_{0,n}(Y,d)]^{\text{virt}}=\EE_d \cap [\overline
M_{0,n}(\PP^s,d)].$
\end{lemma}

At this point we introduce a new object. For any smooth projective
variety $Y$ and any ring $\mathcal A$, we define the formal
completion of $\mathcal A$ along the semigroup of the Mori cone of
$Y$ to be
\begin{eqnarray}
\mathcal A[[q^{\beta}]]:=\{\sum_{\beta}a_{\beta}q^{\beta},  &
a_{\beta}\in \mathcal A, & \beta-\text{effective}\}.
\end{eqnarray}
where $\beta\in H_2(Y,\ZZ)$ is effective if it is a positive
linear combination of algebraic curves. This new ring behaves like
a power series since for each $\beta$, the set of $\alpha$ such
that $\alpha$ and $\beta- \alpha$ are both effective is finite.
For example, in the case of $\PP^s$ we obtain the power series
$\mathcal A[[q]]$.

Choose generators $D_1,...,D_r$ of $H^2(Y,\QQ)$ such that
$j^*(D_1)=H$ and $j^*(D_i)=0$ for $i\geq 2$. Elements of
$H^0(Y,\QQ)\oplus H^2(Y,\QQ)$ are of the form
$t_0+tD:=t_0+t_1D_1+...+t_rD_r$. It is shown in \cite{[15]} that
the generator of the quantum $\mathcal D$-module for the pure
Gromov-Witten theory of $Y$ is
\begin{equation}
J_Y=\text{exp}\left(\frac{t_0+tD}{\hbar}\right)\sum_{\beta\in
H_2(Y,\QQ)}q^{\beta}PD^{-1}\left({e_1}_*\left(\frac{[{\overline
M}_{0,1}(Y,\beta)]^{\text{virt}}}{\hbar(\hbar-c)}\right)\right)
\end{equation}
The moduli spaces $\overline M_{0,1}(Y,\beta)$ are empty unless
$\beta$ is effective. Hence we will consider $J_Y$ as an element
of the ring $H^*Y[[t_0,t_1,...,t_r]][[q^{\beta}]]$.

We extend the map $j^*:H^*Y\rightarrow H^*{\PP^s}$ to a
homomorphism
\begin{equation}
j^*:H^*Y[[t_0,t_1,...,t_r]][[q^{\beta}]]\rightarrow
H^*{\PP^s}[[t_0,t_1]][[q]] \label{map: extend}
\end{equation} by defining $j^*(t_i)=0$ for $i>1$ and $j^*(q^{\beta})=q^{\beta}$ for
$\beta\in j_*(H_2(\PP^s,\ZZ))$ and $j^*(q^{\beta})=0$ for
$\beta\in H_2(Y,\ZZ)-j_*(H_2(\PP^s,\ZZ))$. The following results
show that $J$-function is local.
\begin{theorem}
$j^*(J_Y)=J_V$.
\end{theorem}
{\bf Proof.} Notice that
\begin{equation}
j^*J_Y=\text{exp}\left(\frac{t_0+t_1H}{\hbar}\right)\sum_{d_1=0}^{\infty}q_1^{d_1}j^*PD^{-1}\left({e_1}_*\left(\frac{[{\overline
M}_{0,1}(Y,d_1)]^{\text{virt}}}{\hbar(\hbar-c)}\right)\right)
\end{equation}
Consider the following fiber diagram
\[
\begin{CD}
{\overline M}_{0,1}(\PP^s,d_1) @>{=}>> {\overline M}_{0,1}(Y,d_1) \\
@VVe_1V                                         @VVe_1V\\
\PP^s                 @>{j}>>              Y\\
\end{CD}
\]
By excess intersection theory \cite{[12]} and the previous lemma
\[j^*\left({e_1}_*\left(\frac{[{\overline M}_{0,1}(Y,d_1)]^{\text{virt}}}{\hbar(\hbar-c)}\right)\right)=\EE^-\cap {e_1}_*\left(\frac{\EE_{d_1}}{\hbar(\hbar-c)}\cap [{\overline M}_{0,1}(\PP^s,d_1)]\right).\]
The theorem follows easily.$\dagger$
\begin{theorem}\label{thm: localthm}
Let $V=V^+\oplus V^-=\mathcal O(k)\oplus \mathcal O(-l)$ on
$\PP^s$. Let $\iota:X\hookrightarrow \PP^s$ be the zero locus of a
generic section of $V^+$. Assume that $X$ is smooth and \text{dim}
$X>2$. Let $Y$ be a smooth projective variety such that
$j:X\hookrightarrow Y$ with $\mathcal N={\mathcal
N}_{X/Y}=\iota^*(V^-)$. Assume that if $C\subset Y$ is a curve
with $[C]\in MX$ then all the irreducible components $C_i$ of $C$
satisfy $C_i\subset X$. Let $j^*$ be the map constructed as in
(\ref{map: extend}). Let $J_Y$ be the generator of the pure
$\mathcal D$-module of $Y$ (\cite{[14]}). Then
\[\iota_!(j^*(J_Y))=\text{E}(V^+)J_V.\] where $\iota_!$ is the
Gysin map on cohomology.
\end{theorem}
{\bf Proof.} Since dim$X>2$ it follows that $H^2X$ is generated by
$\iota^*(H)$. Let $\beta_1$ be the Poincar\'{e} dual to
$\iota^*(H)$ and let $D_1,D_2,...,D_r$ be a set of generators of
$H^2(Y,\QQ)$. We may assume that $j^*(D_1)=\iota^*(H)$ and
$j^*(D_i)=0$ for $i>1$. Let $tD:=t_1D_1+...+t_rD_r$. Now
\begin{equation}
J_Y=\text{exp}(\frac{t_0+tD}{\hbar})\sum_{\beta\in
H^2(Y,\QQ)}q^{\beta}PD^{-1}{e_1}_*\left(\frac{[\overline
M_{0,1}(Y,\beta)]^{\text{virt}}}{\hbar(\hbar-c)}\right)
\end{equation}
Consider the following diagram
\[\begin{array}{ccccc}
{\overline M}_{0,1}(Y,d_1\beta_1) & \stackrel{j_1}{\leftarrow} & {\overline M}_{0,1}(X,d_1\beta_1) & \stackrel{\iota_1}{\rightarrow} & {\overline M}_{0,1}(\PP^s,d_1) \\
\downarrow e_1 &                 & \downarrow e_1 &                  & \downarrow e_1 \\
Y &\stackrel{j}{\leftarrow}&  X & \stackrel{\iota}{\rightarrow} &
\PP^s
\end{array}\] The square on the left is a fibre diagram.
We repeatedly use the projection formula:
\[\iota_*(j^*(J_Y))=\iota_!\left(\text{exp}\left(\frac{t_0+t_1\iota^*(H)}{\hbar}\right)\sum_{d_1=0}^{\infty}q_1^{d_1}PD^{-1}j^*{e_1}_*\left(\frac{[\overline M_{0,1}(Y,d_1j_*(\beta_1))]^{\text{virt}}}{\hbar(\hbar-c)}\right)\right)\]
\[=\iota_!\left(\text{exp}\left(\frac{t_0+t_1\iota^*(H)}{\hbar}\right)\sum_{d_1=0}^{\infty}q_1^{d_1}\iota^*(\EE^-)\cup PD^{-1}{e_1}_*{j_1}^*\left(\frac{[\overline M_{0,1}(Y,d_1j_*(\beta_1))]^{\text{virt}}}{\hbar(\hbar-c)}\right)\right)\]
\[=\text{exp}\left(\frac{t_0+t_1H}{\hbar}\right)\left(\EE^++\sum_{d_1=1}^{\infty}q_1^{d_1}(\EE^-)\cup PD^{-1}\iota_*{e_1}_*{j_1}^*\left(\frac{[\overline M_{0,1}(Y,d_1j_*(\beta_1))]^{\text{virt}}}{\hbar(\hbar-c)}\right)\right)\]
\begin{equation}
=\text{exp}\left(\frac{t_0+t_1H}{\hbar}\right)\left(\EE^++\sum_{d_1=1}^{\infty}q_1^{d_1}(\EE^-)\cup
PD^{-1}{e_1}_*{\iota_1}_*{j_1}^*\left(\frac{[\overline
M_{0,1}(Y,d_1j_*(\beta_1))]^{\text{virt}}}{\hbar(\hbar-c)}\right)\right).
\label{eq: J}
\end{equation}
The equality in the second row follows from excess intersection
theory in the left square. An argument similar to
Lemma~\ref{lemma: eqmod} implies that
\begin{equation}
{\overline M}_{0,1}(X,d_1\beta_1)={\overline
M}_{0,1}(Y,d_1j_*(\beta_1)).
\end{equation}
There are two obstruction theories in this moduli stack
corresponding to the moduli problems of maps to $X$ and $Y$
respectively. They differ exactly by the bundle
$R^1{\pi_2}_*e_2^*(\mathcal N)$ where
\[\pi_2:{\overline M}_{0,2}(X,d_1\beta_1)\rightarrow {\overline M}_{0,1}(X,d_1\beta_1)\] is the map that
forgets the second marked point and $\mathcal N=\mathcal N_{X/Y}$.
It follows that:
\[{j_1}^*([\overline
M_{0,1}(Y,d_1j_*(\beta_1))]^{\text{virt}})=\text{E}(R^1{\pi_2}_*e_2^*(\mathcal
N))\cap [\overline M_{0,1}(X,d_1\beta_1)]^{\text{virt}}.\]
Consider the following commutative diagram:
\[
\begin{CD}
{\overline M}_{0,2}(X,d_1\beta_1)   @>{e_2}>>            X \\
@VV{\iota_2}V                                             @VV{\iota}V \\
{\overline M}_{0,2}(\PP^s,d_1)  @>{e_2}>>   \PP^s\\
\end{CD}
\]
We compute:
\[e_2^*(\mathcal N)=e_2^*(\iota^*({\mathcal O}(-l)))={\iota_2}^*{e_2}^*({\mathcal O}(-l)).\]
There is the following fibre square:
\[
\begin{CD}
{\overline M}_{0,2}(X,d_1\beta_1) @>{\iota_2}>> {\overline M}_{0,2}(\PP^s,d_1) \\
@VV{\pi_2}V                                         @VV{\pi_2}V \\
{\overline M}_{0,1}(X,d_1\beta_1)    @>{\iota_1}>> {\overline M}_{0,1}(\PP^s,d_1) \\
\end{CD}
\]
We apply Proposition $9.3$ in \cite{[18]} to obtain:
\[R^1{\pi_2}_*e_2^*(\mathcal N)=R^1{\pi_2}_*{\iota_2}^*\tilde{e_2}^*({\mathcal O}(-l))={\iota_1}^*(R^1{\pi_2}_*\tilde{e_2}^*({\mathcal O}(-l)))={\iota_1}^*(V^-_{d_1}).\]
Therefore:
\begin{eqnarray}
& & {j_1}^*([\overline
M_{0,1}(Y,d_1j_*(\beta_1))]^{\text{virt}})=\text{E}(R^1{\pi_2}_*e_2^*(\mathcal
N))\cap [\overline M_{0,1}(X,d_1\beta_1)]^{\text{virt}}=\nonumber
\\ & & \iota_1^*(E^-_d)\cap [\overline
M_{0,1}(X,d_1\beta_1)]^{\text{virt}}. \label{eq: -J}
\end{eqnarray}
On the other hand, Proposition $11.2.3$ of \cite{[9]} says that :
\begin{eqnarray}
{\iota_1}_*{[\overline
M_{0,1}(X,d_1\beta_1)]^{\text{virt}}}=\EE_d^+\cap [\overline
M_{0,1}(\PP^s,d_1)]. \label{eq: +J}
\end{eqnarray}
Substituting (\ref{eq: -J}) and (\ref{eq: +J}) in (\ref{eq: J}) we obtain
\begin{eqnarray}
& &
\iota_*(j^*(J_Y))=\text{exp}\left(\frac{t_0+t_1H}{\hbar}\right)\cdot
\nonumber \\
& & \cdot
\left(\EE^++\sum_{d=1}^{\infty}q_1^dPD^{-1}{e_1}_*\left(\frac{\EE_d}{\hbar(\hbar-c)}\cap
[\overline M_{0,1}(\PP^s,d)] \right)\cup (\EE^-)\right).
\label{eq: -fin}
\end{eqnarray}

Recall that on $H^*({\overline M}_{0,1}(\PP^s,d))$ we have
$\EE_d=\EE'_d\EE^-_de_1^*(\EE^+)$. Substituting this in (\ref{eq:
-fin}) and using the projection formula we obtain
\begin{eqnarray}
& &
\iota_!(j^*(J_Y))=\text{exp}\left(\frac{t_0+t_1H}{\hbar}\right)\cup
\EE^+\cup \nonumber \\
& & \cup
\left(1+\sum_{d=1}^{\infty}q_1^dPD^{-1}{e_1}_*\left(\frac{\EE'_d\EE^-_d}{\hbar(\hbar-c)}\cap
[\overline M_{0,1}(\PP^s,d)] \right)\cup (\EE^-)\right).
\label{eq: Jbehav}
\end{eqnarray}
The theorem is proven.$\dagger$
\begin{remark}
This naturally leads to local mirror symmetry. For example, let
$Y$ be a Calabi-Yau threefold that contains $X=\PP^2$. By
adjunction formula, the normal bundle of $\PP^2$ in $X$ is
$K_{\PP^2}=\mathcal O_{\PP^2}(-3)$. The last theorem asserts that
the restriction of $J_Y$ in $X$ depends only on $V=\mathcal
O_{\PP^2}(-3)$ i.e. in a neighborhood of $X$ in $Y$. Hence $J_V$
encodes Gromov-Witten correlators of the total space of $\mathcal
O_{\PP^2}(-3)$ which is a local Calabi-Yau. In the next section we
will see that mirror symmetry can be applied to $J_V$ establishing
that mirror symmetry is local at least on the A-side. Interesting
calculations in this direction can be found in \cite{[8]}.
\end{remark}

\section{Mirror Theorem}

In this section we will formulate and prove The Mirror Theorem
which computes the generator $J_V$. Recall that $V=(\oplus_{i\in
I}{\mathcal O(k_i)})\oplus (\oplus_{j\in J}{\mathcal
O(-l_j)})=V^+\oplus V^-$ with $k_i,l_j>0$ for all $i\in I$ and
$j\in J$. Consider the $H^*\PP^s$-valued hypergeometric series
\begin{equation}
I_V(t_0,t_1):=\text{exp}\left(\frac{t_0+t_1H}{\hbar}\right)\sum_{d=0}^{\infty}q^d\frac{\prod_{i\in I} \prod_{m=1}^{k_id}(kH+m\hbar)\prod_{j\in
J}\prod_{m=0}^{l_jd-1}(-l_jH-m\hbar)}{\prod_{m=1}^{d}(H+m\hbar)^{s+1}}. \label{eq: I}
\end{equation}
\begin{theorem}\label{thm: mthm}
({\bf Mirror theorem}). Assume that $\sum_{i\in I}k_i+\sum_{j\in
J}l_j\leq s+1$ and that $J$ is nonempty. If $|J|>1$ or $\sum_{i\in
I}k_i+\sum_{j\in J}l_j<s+1$ then $J_V=I_V$. Otherwise, there
exists a power series $I_1$ of $q$ such that
$J_V(t_0,t_1+I_1)=I_V(t_0,t_1)$ as power series of $q$.
\end{theorem}
\begin{remark}
The case in which $J$ is empty has been treated in \cite {[5]},
\cite{[6]}, \cite {[15]}, \cite{[29]}. It was suggested by
Givental that his techniques should apply in the case in which $J$
is nonempty.
\end{remark}

\subsection{The Equivariant Mirror Theorem.} We use Givental's approach for complete intersections in projective
spaces \cite{[15]} to prove an equivariant version of the theorem.
For the remainder of this paper we will use the standard diagonal
action of $T=(\CC^*)^{s+1}$ on $\PP^s$ with weights
$(-\lambda_0,...,-\lambda_s)$. Recall from section $2.2$ that
$\mathcal
P=H^*_T(\PP^s)={\CC}[\lambda][p]/\prod_{i=0}^{s}(p-{\lambda}_i)$
and $\mathcal R=\CC(\lambda)[p]/\prod_{i=0}^{s}(p-{\lambda}_i)$.
Denote
\begin{eqnarray}
J_V^{eq}:=\text{exp}\left(\frac{t_0+pt_1}{\hbar}\right)\sum_{d=0}^{\infty}q^d{e_1}_*\left(\frac{E'_dE_d^-
}{\hbar(\hbar-c)}\right)\cup (\prod_{j\in J}-l_jp)= \nonumber \\
\text{exp}\left(\frac{t_0+t_1p}{\hbar}\right)S(q,\hbar).
\label{eq: Jeq}
\end{eqnarray}
and
\begin{eqnarray}
I_V^{eq}:=\text{exp}\left(\frac{t_0+t_1p}{\hbar}\right)\sum_{d=0}^{\infty}q^d\frac{\prod_{i\in
I}\prod_{m=1}^{k_id}(k_ip+m\hbar)\prod_{j\in
J}\prod_{m=0}^{l_jd-1}(-l_jp-m\hbar)}{\prod_{m=1}^{d}\prod_{i=0}^{s}(p-\lambda_i+m\hbar)}=
\nonumber \\
\text{exp}\left(\frac{t_0+t_1p}{\hbar}\right)S'(q,\hbar).
\label{eq: Ieq}
\end{eqnarray}

Obviously the nonequivariant limits of $J_V^{eq}$ and $I_V^{eq}$
are respectively $J_V$ and $I_V$. The mirror theorem will follow
as a nonequivariant limit of the following
\begin{theorem}
({\bf The Equivariant Mirror Theorem}). The same change of
variables from theorem $5.0.3$ transforms $I^{eq}_V$ into
$J^{eq}_V$.
\end{theorem}
\begin{remark}
As the reader will see, the central part of the proof of the
Mirror Theorem (up to section $5.5$) involves lengthy formulaes
and algebraic manipulations. To simplify the presentation, we will
assume during this part that $V=\mathcal O(k)\oplus \mathcal
O(-l)$. The general case is similar. We will return to the general
case $V=V^+\oplus V^-$ in section $5.5$.
\end{remark}
Recall that the equivariant Thom classes $\phi_i$ of the fixed
points $p_i$ form a basis of $\mathcal R$ as a
$\CC(\lambda)$-vector space. Let $S_i$ and  $S'_i$ be the
restrictions of $S$ and $S'$ at the fixed point $p_i$. By the
localization theorem in $\PP^s$ they determine $S$ and $S'$. By
the projection formula
\begin{eqnarray}
S_i=\int_{\PP^s_T}S\cup \phi_i=1+\sum_{d=1}^{\infty}q^d
\int_{{\overline M}_{0,1} (\PP^s,d)_T}
\frac{e_1^*(-lp\phi_i)}{\hbar(\hbar-c)}E'_dE^-_d. \label{eq:
Siforlater}
\end{eqnarray}
The proof of the equivariant mirror theorem is based on exhibiting
similar properties of the correlators $S_i$ and $S'_i$. The extra
property $S_i=1+o(\hbar^{-2})$ determines $S_i$ uniquely. After
the change of variables that property is satisfied by $S'_i$ as
well which implies $S_i=S'_i$.

We now proceed with displaying properties of the correlators $S_i$ and $S'_i$.

\subsection{\bf Linear recursion relations} The first property is given by this
\begin{lemma} \label{lemma: recursionlemma}
The correlators $S_i$ satisfy the following linear recursion relations:
\begin{equation}
S_i=1+\sum_{d=1}^{\infty}q^dR_{id}+\sum_{d=1}^{\infty}\sum_{j\neq i}q^dC_{ijd}S_j\left(q,\frac{\lambda_j-\lambda_i}{d}\right)
\end{equation} where $R_{id}\in \CC(\lambda)[\hbar^{-1}]$ are polynomials in $\hbar^{-1}$ and
\begin{equation}
C_{ijd}=\frac{(\lambda_j-\lambda_i)\prod_{m=1}^{kd}(k\lambda_i+m\frac{\lambda_i-\lambda_j}{d} )\prod_{m=0}^{ld-1}(-l\lambda_i+m\frac{\lambda_i-\lambda_j}{d})}{d\hbar(d\hbar+\lambda_i-\lambda_j)\prod_{m=1}^{d}\prod_{k=0,(k,m)\neq (j,d)}^{s}(\lambda_i-\lambda_k+m\frac{\lambda_j-\lambda_i}{d})}.
\end{equation}
\end{lemma}
{\bf Proof.} We will see during the proof that $S_j$ is regular at $\displaystyle{\hbar=\frac{\lambda_j-\lambda_i}{d}}$.
The integrals that appear in the formula for $S_i$ can be evaluated using localization theorem
\begin{equation}
\int_{{\overline M}_{0,1}
({\PP}^s,d)_T}\frac{e_1^*(-lp\phi_i)}{\hbar(\hbar-c)}E'_dE^-_d=\sum_{\Gamma}\int_{(\overline{M}_{\Gamma})_T}\frac{1}{a_{\Gamma}\text{E}_T(N_{\Gamma})}
\left(\frac{e_1^*(-lp\phi_i)}{\hbar(\hbar-c)}E'_dE^-_d\right)_{\Gamma}
\end{equation}
There are three types of fixed point components $M_{\Gamma}$ of
${\overline M}^T_{0,1} ({\PP}^s,d)$. The first one consists of
those $M_{\Gamma}$ where the component of the curve that contains
the marked point is collapsed to $p_i$. We denote the set of these
components by $\bf {\mathcal F}^i_{1,d}$. Let $\bf {\mathcal
F}^i_{2,d}$ be the set of those  $M_{\Gamma}$ in which the marked
point is mapped at $p_i$ and its incident component is a multiple
cover of the line $\overline{p_i,p_j}$ for some $j\neq i$. Finally
let $\bf {\mathcal F}^i_{0,d}$ be the rest of the fixed point
components. Notice first that:
\begin{equation}
\sum_{\Gamma\in \ {\mathcal F}^i_{0,d}}\int_{({\overline
M}_{\Gamma})_T}{\frac{1}{a_{\Gamma}\text{E}_T(N_{\Gamma})}
\left(\frac{e_1^*(-lp\phi_i)}{\hbar(\hbar-c)}E'_dE^-_d\right)_{\Gamma}}=0.
\label{eq: nocont}
\end{equation}
Indeed, let $\Gamma_j\in \bf {\mathcal F}^i_{0,d}$ represent a
fixed point component with the marked point mapped to the fixed
point $p_j$ for some $j\neq i$. Since
$(e_1^*(\phi_i))_{\Gamma_j}=0$ we are done. Next, in each fixed
point component that belongs to $\bf {\mathcal F}^i_{1,d}$ the
class $c$ is nilpotent. Indeed, if $\Gamma$ is the decorated graph
that represents such a fixed point component, let ${\overline
M}_{0,k}$ correspond to the vertex of $\Gamma$ that contains the
marked point. Then $k\leq d+1$. There is a morphism:
\[\varphi:M_{\Gamma}\mapsto {\overline M}_{0,k}\] such that $\varphi^*(c)=c_{\Gamma}$. For dimension reasons
$c^{d-1}=0$ on ${\overline M}_{0,k}$ therefore $\displaystyle{\frac{1}{\hbar(\hbar-c)}}$ is a polynomial of $c$ in ${\overline M}_{\Gamma}$.
Hence
\begin{equation}
\sum_{\Gamma\in \ {\mathcal
F}^i_{1,d}}\int_{(\overline{M}_{\Gamma})_T}{\frac{1}{a_{\Gamma}\text{E}_T(N_{\Gamma})}
\left(\frac{e_1^*(-lp\phi_i)}{\hbar(\hbar-c)}E'_dE^-_d\right)_{\Gamma}}=R_{id}
\label{eq: polcont}  \end{equation} is a polynomial in
$\hbar^{-1}$.

We now consider the fixed point components in $\bf {\mathcal
F}^i_{2,d}$. Again let $\Gamma$ represent such a component. For a
stable map $(C,x_1,f)$ in $\Gamma$ let $C'$ be the component of
$C$ containing $x_1$, $C''$ the rest of the curve, $x=C'\cap C''$
and $f(x)=p_j$ for some $j\neq i$. Let $d'$ be the degree of the
map f on the component $C'$ and $d''=d-d'$. Then
$(C'',x,f|_{C''})$ is a fixed point in   ${\overline M}_{0,1}
({\PP}^s,d'')$. Denote its decorated graph by ${\Gamma}''$. Choose
the coordinates on $C'$ such that the restriction of $f$ on $C'$
is given by
$f(y_0,y_1)=(0,...,z_i=y_0^{d'},...,z_j=y_1^{d'},...,0)$. As
$\Gamma$ moves in $\bf {\mathcal F}^i_{2,d}$, the set of all such
${\Gamma}''$ exhausts all the fixed points in ${\overline M}_{0,1}
({\PP}^s,d'')$ where the first marked point is not mapped to
$p_i$. Since $\text{Aut}(\Gamma)=\text{Aut}(\Gamma'')$ it follows
from (\ref{eq: autorder}) that:
\begin{eqnarray}
a_{\Gamma}=d'a_{{\Gamma}''}.
\end{eqnarray}
The local coordinate at $p_i$ on the component $C'$ is
$z=\displaystyle{\frac{y_1}{y_0}}$. The weight of the $T$-action
on $y_l$ is $\displaystyle \frac{\lambda_l}{d'}$ for $l=0,1$. It
follows that the weight of the action on the coordinate $z$ and
hence on $T^*_{p_i}C'$ is
$\displaystyle{\frac{\lambda_j-\lambda_i}{d'}}$ therefore
$\displaystyle{c_{\Gamma}=\frac{\lambda_j-\lambda_i}{d'}}$. Now
$\text{E}_T(N_{\Gamma})$ can be split in three pieces: smoothing
the node $x$ and deforming the maps $f|_{C''}$ and $f|_{C'}$. It
follows \cite{[15]}
\begin{eqnarray}
& &
\text{E}_T(N_{\Gamma})=\left(\frac{\lambda_j-\lambda_i}{d'}-c_{{\Gamma}}''\right)\text{E}_T(N_{{\Gamma}''})\cdot
\nonumber \\
& & \cdot \prod_{m=0}^{d'-1}\prod_{k=0,(m,k)\neq
(0,i)}^{s}\left(\lambda_i-\lambda_k+m\frac{\lambda_j-\lambda_i}{d'}\right).
\label{eq: ENormalnode}
\end{eqnarray}
Next, we find the localization of $E'_d$ and $E^-_d$ on the fixed
point component $M^i_{d_1d_2}$. Consider the normalization
sequence at the node $x$
\begin{equation}
0\rightarrow {\mathcal O}_C\rightarrow {\mathcal O}_{C'}\oplus {\mathcal O}_{C''}\rightarrow {\mathcal O}_x\rightarrow 0. \label{eq: normal}
\end{equation}
Twisting it by $f^*(V^+)$ and $f^*(V^-)$ respectively  and taking the cohomology sequence yields
\[(E^-_d)_{\Gamma}=(-l\lambda_j)(E^-_{d''})_{{\Gamma}''}(E^-_{d'})_{{\Gamma}'}\] and
\begin{equation}
(E^+_d)_{\Gamma}=\frac{(E^+_{d'})_{{\Gamma}'}(E^+_{d''})_{{\Gamma}''}}{k\lambda_j}
\label{eq: splitE+}
\end{equation}
An explicit basis for $H^1(C',f^*(V^-))=H^1({\mathcal
O}_{{\PP}^1}(-ld')$ consists of
\begin{eqnarray}
\frac{y_0^sy_1^{ld'-2-s}}{(y_0y_1)^{ld'-1}}=\frac{1}{y_0^{ld'-s-1}y_1^{1+s}}
: s=0,1,...,ld'-2.
\end{eqnarray}
It allows us to compute:
\[(E^-_{d'})_{{\Gamma}'}=\prod_{s=0}^{ld'-2}\left(\frac{1+s-ld'}{d'}\lambda_i-\frac{1+s}{d'}\lambda_j\right)=\prod_{s=1}^{ld'-1}\left(-l\lambda_i+s\frac{\lambda_i-\lambda_j}{d'}\right).\]
Therefore we have:
\begin{equation}
(E^-_d)_{\Gamma}=(-l\lambda_j)\prod_{s=1}^{ld'-1}\left(-l\lambda_i+s\frac{\lambda_i-\lambda_j}{d'}\right)(E^-_{d''})_{{\Gamma}''}. \label{eq: E-}
\end{equation}
A basis for $H^0(C',f^*(V^+))=H^0({\mathcal O}_{{\PP}^1}(kd'))$
consists of monomials $y_0^sy_1^{kd'-s}$ for $s=0,...kd'$. It can
be used to calculate $(E^+_{d'})_{{\Gamma}'}$ similarly to
$(E^-_{d'})_{{\Gamma}'}$ above. Recall from (\ref{eq: Eulersplit})
that $E^+_d=e_1^*(\text{E}(V^+))E'_d$. The line bundle
$e_1^*(V^+)$ is trivial on $M_{\Gamma}$ but the torus acts on it
with weight $k\lambda_i$. Hence
$(E^+_d)_{\Gamma}=k\lambda_i(E'_d)_{\Gamma}$. Substituting in
(\ref{eq: splitE+}) yields
\begin{equation}
(E'_d)_{\Gamma}=\prod_{r=1}^{kd'}
\left(k\lambda_i+r\frac{\lambda_j-\lambda_i}{d'}\right)(E'_{d''})_{{\Gamma}''}. \label{eq: E'd}
\end{equation}
We pause here to show that $S_i$ is regular at $\displaystyle{\hbar=\frac{\lambda_j-\lambda_i}{d}}$ for any $j\neq i$ and any $d>0$.
It follows from (\ref{eq: nocont}) and (\ref{eq: polcont}) that:
\begin{equation}
S_i=1+\sum_{d=1}^{\infty}q^dR_{id}+\sum_{\Gamma\in \ {\mathcal
F}^i_{2,d}}q^d
\int_{(\overline{M}_{\Gamma})_T}\frac{(-l\lambda_i)\prod_{k\neq
i}(\lambda_i-\lambda_k)(E'_d\cdot
E^-_d)_{\Gamma}}{\hbar(\hbar-c_{\Gamma})a_{\Gamma}\text{E}_T(N_{\Gamma})}.
\label{eq: S_i}
\end{equation}
From this representation of $S_i$ it is clear that the
coefficients of the power series
$S_i=\sum_{d=0}^{\infty}S_{id}q^d$ belong to $\QQ(\lambda,\hbar)$.
But $\displaystyle{c_{\Gamma}=\frac{\lambda_j-\lambda_i}{d'}}$ for
some $d'\leq d$ and $R_{id}$ has poles only at $\hbar=0$ therefore
$S_i$ is regular at
$\displaystyle{\hbar=\frac{\lambda_i-\lambda_j}{d}}$. We use the
equations (\ref{eq: ENormalnode}), (\ref{eq: E'd}) and (\ref{eq:
E-}) to compute
\begin{eqnarray}
& & \sum_{\Gamma\in \ {\mathcal F}^i_{2,d}}q^d\int_{(\overline{M}_{\Gamma})_T} \frac{(-l\lambda_j)\prod_{k\neq i}(\lambda_i-\lambda_k)(E'_d\cdot E^-_d)_{\Gamma}}{\hbar(\hbar-c_{\Gamma})a_{\Gamma}\text{E}_T(N_{\Gamma})}= \nonumber \\
& & \sum_{d'=1}^{\infty}\sum_{j\neq i}q^{d'}C_{ijd'}\sum_{\Gamma''}q^{d''}\int_{(\overline{M}_{{\Gamma}''})_T}\frac{-l\lambda_j\prod_{k\neq j}(\lambda_j-\lambda_k)(E'_{d''}\cdot E^-_{d''})_{{\Gamma}''}}{(\frac{\lambda_j-\lambda_i}{d'})(\frac{\lambda_j-\lambda_i}{d'}-c_{\Gamma}'')a_{{\Gamma}''}\text{E}_T(N_{{\Gamma}''})}= \nonumber \\
& & \sum_{d'=1}^{\infty}\sum_{j\neq i}q^{d'}C_{ijd'}S_j\left(q,\frac{\lambda_j-\lambda_i}{d'}\right).
\end{eqnarray}
The lemma follows by substituting the above identity into (\ref{eq: S_i})$\dagger$.
\smallskip
\begin{lemma}
The correlators $S'_i$ satisfy the same linear recursion relations
as $S_i$.
\end{lemma}
{\bf Proof.} We know that  $S'_i=\sum_{d=0}^{\infty}q^dS'_{id}$ with
\begin{equation}
S'_{id}=\frac{\prod_{m=1}^{kd}(k\lambda_i+m\hbar)\prod_{m=0}^{ld-1}(-l\lambda_i-m\hbar)}{d\hbar\prod_{m=1}^{d}\prod_{j=0,(j,m)\neq (i,d)}^{s}(\lambda_i-{\lambda}_j+m\hbar)}  \label{eq: S'_i}
\end{equation}
Note that $S'_{id}\in \CC(\lambda,\hbar)$ is a proper rational expression of $\hbar$. It has multiple poles
at $\hbar=0$ and simple poles at $\hbar=\displaystyle{\frac{\lambda_r-\lambda_i}{m}}$ for any $r\neq i$ and
any $1\leq m\leq d$. Applying calculus of residues in the $\hbar$-variable yields:
\begin{eqnarray}
& & S'_{id}=R_{id}+\sum_{m=1}^{d}\sum_{r\neq i}\frac{1}{d\hbar(\lambda_i-\lambda_r+m\hbar)} \nonumber \\ & & \times \frac{\prod_{n=1}^{kd}(k\lambda_i+n\frac{\lambda_r-{\lambda}_i}{m})\prod_{n=0}^{ld-1}(-l\lambda_i-n\frac{\lambda_r-{\lambda}_i}{m})}{\prod_{n=1,(j,n)\neq (r,m)}^{d}\prod_{j=0,(j,n)\neq (i,d)}^{s}(\lambda_i-{\lambda}_j+n  \frac{\lambda_r-{\lambda}_i}{m})} \label{eq: res}
\end{eqnarray}
for some polynomials $R_{id}\in \CC(\lambda)[\hbar^{-1}]$ such that $R_{id}(0)=0$. Substitute equation (\ref{eq: res}) in (\ref{eq: S'_i}) to obtain:
\begin{eqnarray}
& & S'_i=1+\sum_{d=1}^{\infty}q^dR_{id}+\sum_{d=1}^{\infty}q^d\sum_{r\neq i}\sum_{m=1}^{d}\frac{1}{d\hbar(\lambda_i-\lambda_r+m\hbar)} \nonumber \\ & & \times \frac{\prod_{n=1}^{kd}(k\lambda_i+n\frac{\lambda_r-{\lambda}_i}{m})\prod_{n=0}^{ld-1}(-l\lambda_i-n\frac{\lambda_r-{\lambda}_i}{m})}{\prod_{n=1,(j,n)\neq (r,m)}^{d}\prod_{j=0,(j,n)\neq (i,d)}^{s}(\lambda_i-{\lambda}_j+n\frac{\lambda_r-{\lambda}_i}{m})}.
\end{eqnarray}
Changing the order of summation in the last equation yields:
\begin{eqnarray}
& & S'_i-1-\sum_{d=1}^{\infty}q^dR_{id}=\sum_{r\neq i}\sum_{m=1}^{\infty}q^m\frac{1}{\hbar(\lambda_i-\lambda_r+m\hbar)} \nonumber \\ & & \times \sum_{d=m}^{\infty}q^{d-m}\frac{\prod_{n=1}^{kd}(k\lambda_i+n\frac{\lambda_r-{\lambda}_i}{m}    )\prod_{n=0}^{ld-1}(-l\lambda_i-n\frac{\lambda_r-{\lambda}_i}{m}    )}{d\prod_{n=1,(j,n)\neq (r,m)}^{d}\prod_{j=0,(j,n)\neq (i,d)}^{s}(\lambda_i-{\lambda}_j+n  \frac{\lambda_r-{\lambda}_i}{m})}.
\end{eqnarray}
The lemma follows from the identity
\begin{eqnarray}
& &
\sum_{d=m}^{\infty}q^{d-m}\frac{\prod_{n=1}^{kd}(k\lambda_i+n\frac{\lambda_r-{\lambda}_i}{m})\prod_{n=0}^{ld-1}(-l\lambda_i-n\frac{\lambda_r-{\lambda}_i}{m}
)}{d\prod_{n=1,(j,n)\neq (r,m)}^{d}\prod_{j=0,(j,n)\neq
(i,d)}^{s}(\lambda_i-{\lambda}_j+n
\frac{\lambda_r-{\lambda}_i}{m})} \nonumber \\
& &
=\frac{(\frac{\lambda_r-\lambda_i}{m})\prod_{n=1}^{km}(k\lambda_i+n\frac{\lambda_r-{\lambda}_i}{m})\prod_{n=0}^{lm-1}(-l\lambda_i-n\frac{\lambda_r-{\lambda}_i}{m})}{\prod_{n=1,(j,n)\neq
(r,m)}^{m}\prod_{j=0}^{s}(\lambda_i-{\lambda}_j+n
\frac{\lambda_r-{\lambda}_i}{m})} \nonumber \\
& & \times
\sum_{u=0}^{\infty}q^u\frac{\prod_{n=1}^{ku}(k\lambda_r+n\frac{\lambda_r-{\lambda}_i}{m})\prod_{n=0}^{lu-1}(-l\lambda_r-n\frac{\lambda_r-{\lambda}_i}{m})}{u(\frac{\lambda_r-\lambda_i}{m})\prod_{n=1,(j,n)\neq
(r,u)}^{u}\prod_{j=0}^{s}\left(\lambda_r-{\lambda}_j+n
\frac{\lambda_r-{\lambda}_i}{m}\right)}
\end{eqnarray}

\subsection{\bf Double polynomiality}
Recall from Section $3.1$ that $V$ induces a modified equivariant
integral $\omega_V: \mathcal R\rightarrow \CC(\lambda)$ defined as
follows:
\begin{equation}
\omega_V(a):=\int_{\PP^s_{\bf T}}a\cup \frac{E^+}{E^-}.
\end{equation}
As one can see, in the case $V=\mathcal O(k)\oplus \mathcal
O(-l)$, this modified equivariant integral simplifies via
$\displaystyle{\frac{E^+}{E^-}=\frac{kp}{-lp}=\frac{k}{l}}$. We
have chosen not to simplify this integral in the proof of the
following Lemma so that it is easier to see how to proceed\ in the
general case.
\begin{lemma}
If $z$ is a variable, the expression:
\[P(z,\hbar)=\omega_V\left(e^{pz}S(qe^{z\hbar},\hbar)S(q,-\hbar)\right)\]
belongs to $\QQ(\lambda)[\hbar][[q,z]].$
\end{lemma}
{\bf Proof.}  In Section $2.3$ we introduced the action of
$T'=T\times \CC^*$ on $\PP^s\times \PP^1$ with weights
$(-\lambda_0,...,-\lambda_s)$ on the first factor and $(-\hbar,0)$
in the second factor. Consider the following $T'$-equivariant
diagram
\[ \begin{CD}
{\overline M}_{0,1}({\PP}^s\times{\PP}^1,(d,1))@>e_1>>{\PP}^s\times{{\PP}^1} \\
@VV \pi V \\
M_d={\overline M}_{0,0}({\PP}^s\times{\PP}^1,(d,1))
\end{CD} \] Define:
\begin{equation}
W_d=W^+_d\oplus W^-_d:=\pi_*((e_1)^*({\mathcal O}_{\PP^s}(k)\otimes{\mathcal O}_{{\PP}^1}))\oplus R^1\pi_*((e_1)^*({\mathcal O}_{\PP^s}(-l)\otimes{\mathcal O}_{{\PP}^1}))
\end{equation}
The lemma will follow from the identity:
\begin{equation}
P(z,\hbar)=\sum_{d=0}^{\infty}q^d\int_{(M_d)_{T'}}e^{z{\psi}^*{\kappa}}\text{E}_{T'}(W_d)
\label{eq: auxid}
\end{equation}
where $\psi$ and $\kappa$ were defined in Section $2.3$. The
localization formula for the diagonal action of $T$ on $\PP^s$
applied to the left side gives
\begin{equation}
P(z,\hbar)=\sum_{i=0}^{s}\frac{S_i(qe^{z\hbar},\hbar)e^{z\lambda_i}S_i(q,-\hbar)}{\prod_{k\neq i}(\lambda_i-\lambda_k)}\left(\frac{k\lambda_i}{-l\lambda_i}\right). \label{eq: WfromPP}
\end{equation}
We recall from identity (\ref{eq: Siforlater})
\begin{equation}
S_i=1+\sum_{d=1}^{\infty}q^d\int_{{\overline M}_{0,1}
({\PP}^s,d)_T}\frac{e_1^*(-lp\phi_i)}{\hbar(\hbar-c)}E'_dE^-_d.
\label{eq: recSi}
\end{equation}
To compute the integrals on the right side of (\ref{eq: auxid}) we
will use localization for the action of $T'$ on $M_d$. In Section
$2.3$ we found that the components of the fixed point loci have
the form $M^i_{d_1d_2}={\overline M}_{\Gamma^i_{d_1}}\times
\overline{M}_{\Gamma^i_{d_2}}$ for some $i=0,1,...,s$ and a
splitting $d=d_1+d_2$. We first compute the restriction of
$\text{E}_{T'}(W_d)$ in such a component. Consider the following
normalization sequence
\begin{equation}
0\rightarrow \mathcal O_C\rightarrow \mathcal O_{C_0}\oplus \mathcal O_{C_1}\oplus \mathcal O_{C_2}\rightarrow \mathcal O_{x_1}\oplus \mathcal O_{x_2}\rightarrow 0 \label{eq: eqnormalization}
\end{equation}
Twist (\ref{eq: eqnormalization}) by $f^*({\mathcal
O}(-l)\otimes{\mathcal O}_{{\PP}^1})$ and take the corresponding
long exact cohomology sequence. We obtain
\[0\rightarrow {\CC}\rightarrow {\mathcal O}_{x_1}(-l)\oplus {\mathcal O}_{x_2}(-l)\rightarrow W^-_d\rightarrow W^-_{d_1}\oplus W^-_{d_2}\rightarrow 0.\] The first piece is trivial since it comes from the isomorphism
\[({\mathcal O}_{\PP^s}(-l)\times {\mathcal O}_{\PP^1})|_{C_0}\cong {\mathcal O}_{C_0}\cong {\CC}.\]
The left hand side is generated by $\displaystyle{\frac{1}{{z_i}^l}}$ therefore the weight of that piece
is $-l\lambda_i$. It follows that
\[\text{E}_{T'}(W^-_d)=(-l\lambda_i)E^-_{d_1}E^-_{d_2}.\]
Similarly, twisting the normalization sequence (\ref{eq: eqnormalization})
by $f^*({\mathcal O}(k)\otimes{\mathcal O}_{{\PP}^1})$ and taking the corresponding cohomology sequence we obtain:
\[\text{E}_{T'}(W^+_d)=(k\lambda_i)E'_{d_1}E'_{d_2}.\]
We now use the localization theorem to calculate the integrals on
the right side of (\ref{eq: auxid}). The equivariant Euler class
of the normal bundle of the fixed point component $M^i_{d_1d_2}$
has been calculated at the end of section $2.3$.
\begin{eqnarray}
\int_{(M_d)_{T'}}e^{z{\psi}^*{\kappa}}\text{E}_{T'}(W_d)=\sum_{\Gamma^i_{d_1},\Gamma^i_{d_2}}(k\lambda_i)(-l\lambda_i)\frac{e^{z(\lambda_i+d_2\hbar)}}{\prod_{k\neq
i}(\lambda_i-\lambda_k)}  \\
 \times    \int_ {(\overline{M}_{\Gamma^i_{d_1}})_T}
\frac{1}{E_T(N_{\Gamma^i_{d_1}})}
\left(\frac{e_1^*(\phi_i)E'_{d_1}E^-_{d_1}}{-\hbar(-\hbar-c_1)}\right)_{\Gamma^i_{d_1}}
\int_{(\overline{M}_{\Gamma^i_{d_2}})_T}  {
\frac{1}{E_T(N_{\Gamma^i_{d_2}})}
\left(\frac{e_1^*(\phi_i)E'_{d_2}E^-_{d_2}}{\hbar(\hbar-c_2)}\right)_{\Gamma^i_{d_2}}
}= \nonumber \\
\sum_{\Gamma^i_{d_1},\Gamma^i_{d_2}}\frac{k\lambda_i}{-l\lambda_i}\frac{e^{z(\lambda_i+d_2\hbar)}}{\prod_{k\neq
i}(\lambda_i-\lambda_k)} \int_{(\overline{M}_{\Gamma^i_{d_1}})_T}
\frac{1}{E_T(N_{\Gamma^i_{d_1}})}\left(\frac{e_1^*(-l\lambda_i\phi_i)E_{d_1}}{-\hbar(-\hbar-c_1)}
\right)_{\Gamma^i_{d_1}} \times \nonumber \\
\times \int_{(\overline{M}_{\Gamma^i_{d_2}})_T}
\frac{1}{E_T(N_{\Gamma^i_{d_2}})}\left(\frac{e_1^*(-l\lambda_i\phi_i)E_{d_2}}{\hbar(\hbar-c_2)}
\right)_{\Gamma^i_{d_2}}. \nonumber
\end{eqnarray}
If we use localization to compute $S_i$ in (\ref{eq: recSi}) and
then substitute in (\ref{eq: WfromPP}) we obtain the right side of
the last equation.$\dagger$
\begin{lemma} \label{lemma: doublelemma}
If $z$ is a variable, the expression:
\begin{eqnarray}
P'(z,h)=\omega_V\left(S'(qe^{z\hbar},\hbar)e^{pz}S'(q,-\hbar)\right)
\label{eq: doubleS'}
\end{eqnarray}
belongs to $\QQ(\lambda)[\hbar][[q,z]].$
\end{lemma}
{\bf Proof.} The lemma will follow from the identity
\begin{equation}
P'(z,h)=\sum_{d=0}^{\infty}q^d\int_{(N_d)_{T'}}{e^{z\kappa}\prod_{m=0}^{kd}(k\kappa-m\hbar)\prod_{m=1}^{ld-1}(-l\kappa+m\hbar)}.
\label{eq: doubleNd}
\end{equation}
For $d=0$ the convention
\[\int_{(N_d)_{T'}}{e^{z\kappa}\prod_{m=0}^{kd}(k\kappa+m\hbar)\prod_{m=1}^{ld-1}(-l\kappa+m\hbar)}=\int_{\PP^s_{\bf T}}e^{pz}\left(\frac{\prod_{i\in I} k_ip}{\prod_{j\in J} -l_jp}\right)\] is taken.
Apply the localization formula to the integrals (\ref{eq:
doubleS'}).
\begin{eqnarray}
& &
P'(z,\hbar)=\sum_{i=0}^{s}\frac{k\lambda_ie^{\lambda_iz}}{(-l\lambda_i)\prod_{j\neq
i}(\lambda_i-\lambda_j)} \nonumber \\
& & \times
\sum_{d_1=0}^{\infty}(qe^{z\hbar})^{d_1}\frac{\prod_{m=1}^{kd_1}(k\lambda_i+m\hbar)\prod_{m=0}^{ld_1-1}(-l\lambda_i-m\hbar)}{\prod_{m=1}^{d_1}\prod_{j=0}^{s}(\lambda_i-{\lambda}_j+m\hbar)}
\nonumber \\
& & \times
\sum_{d_2=0}^{\infty}q^{d_2}\frac{\prod_{m=1}^{kd_2}(k\lambda_i-m\hbar)\prod_{m=0}^{ld_2-1}(-l\lambda_i+m\hbar)}{\prod_{m=1}^{d_2}\prod_{j=0}^{s}(\lambda_i-{\lambda}_j-m\hbar)}=
\nonumber \\
& & \sum_{i=0}^{s}\frac{1}{\prod_{j\neq
i}(\lambda_i-\lambda_j)}\sum_{d_1=0}^{\infty}q^{d_1z}e^{(\lambda_i+d_1\hbar)z}\frac{\prod_{m=0}^{kd_1}(k\lambda_i+m\hbar)\prod_{m=1}^{ld_1-1}(-l\lambda_i-m\hbar)}{\prod_{m=1}^{d_1}\prod_{j=0}^{s}(\lambda_i-{\lambda}_j+m\hbar)}
\nonumber \\
& & \times
\sum_{d_2=0}^{\infty}q^{d_2}\frac{\prod_{m=1}^{kd_2}(k\lambda_i-m\hbar)\prod_{m=0}^{ld_2-1}(-l\lambda_i+m\hbar)}{\prod_{m=1}^{d_2}\prod_{j=0}^{s}(\lambda_i-\lambda_j-m\hbar)}.
\end{eqnarray}
But for $d_1,d_2>0$
\begin{eqnarray}
& & \frac {\prod_{m=0}^{kd_1}(k\lambda_i+m\hbar)\prod_{m=1}^{ld_1-1}(-l\lambda_i-m\hbar)\prod_{m=1}^{kd_2}(k\lambda_i-m\hbar)\prod_{m=0}^{ld_2-1}(-l\lambda_i+m\hbar)}
{\prod_{j\neq i}(\lambda_i-\lambda_j)\prod_{m=1}^{d_1}\prod_{j=0}^{s}(\lambda_i-{\lambda}_j+m\hbar)\prod_{m=1}^{d_2}\prod_{j=0}^{s}(\lambda_i-\lambda_j-m\hbar)}= \nonumber \\ & & \frac{\prod_{m=0}^{k(d_1+d_2)}(k(\lambda_i+d_1\hbar)-m\hbar)\prod_{m=1}^{l(d_1+d_2)-1}(-l(\lambda_i+d_1\hbar)+m\hbar)}{\prod_{j=0}^{s}\prod_{m=0,(j,m)\neq (i,d_1)}^{d_1+d_2}(\lambda_i+d_1\hbar-\lambda_j-m\hbar)}.\nonumber
\end{eqnarray}
Therefore
\begin{eqnarray}
& & P'(z,\hbar)=\sum_{d=0}^{\infty}q^d\sum_{d_1=0}^{d}\sum_{i=0}^{s}e^{(\lambda_i+d_1\hbar)z} \nonumber \\ & & \times \frac{\prod_{m=0}^{kd}(k(\lambda_i+d_1\hbar)-m\hbar)\prod_{m=1}^{ld-1}(-l(\lambda_i+d_1\hbar)+m\hbar)}{\prod_{j=0}^{s}\prod_{m=0,(j,m)\neq (i,d_1)}^{d}(\lambda_i+d_1\hbar-\lambda_j-m\hbar)}.
\end{eqnarray}
By the localization formula in $N_d$ we can see that
\[P'(z,\hbar)=\sum_{d=0}^{\infty}q^d\int_{(N_d)_{T'}}{e^{z\kappa}\prod_{m=0}^{kd}(k\kappa-m\hbar)\prod_{m=1}^{ld-1}(-l\kappa+m\hbar)}.\]

The lemma is proven.$\dagger$

\subsection{\bf Mirror transformation and uniqueness} The following two theorems carry over from \cite{[15]}.
The first lemma deals with uniqueness.
\begin{lemma}
Let $S=\sum_{d=0}^{\infty}S_dq^d$ and
$S'=\sum_{d=0}^{\infty}S'_dq^d$ be two power series with
coefficients in ${\mathcal R}[[\hbar^{-1}]]$ that satisfy the
following conditions:
\begin{enumerate}
\item $S_0=S'_0=1$
\item They both satisfy the recursion relations of Lemma \ref{lemma: recursionlemma}.
\item They both have the double polynomiality property of Lemma \ref{lemma: doublelemma}.
\item For any d, $S_d\equiv S'_d$ mod $(\hbar^{-2})$.
\end{enumerate}
Then $S=S'$.
\end{lemma}

The second lemma describes a transformation which preserves the
properties of lemma $5.4.1$.
\begin{lemma}
Let $I_1$ be a power series in $q$ whose first term is zero. Then
$\text{exp}(\frac{I_1p}{\hbar})S(qe^{I_1},\hbar)$ satisfies
conditions $1,2,3$ of Lemma $5.4.1$.
\end{lemma}

\subsection {The conclusion of the proof of Mirror Theorem} Recall that
\[I_V^{eq}=\text{exp}\left(\frac{t_0+t_1p}{\hbar}\right)\left(1+\sum_{d=1}^{\infty}q^d\frac{\prod_{i\in I}\prod_{m=1}^{k_id}(k_ip+m\hbar)\prod_{j\in J}\prod_{m=0}^{l_jd-1}(-l_jp-m\hbar)}{\prod_{m=1}^{d}\prod_{i=0}^{s}(p-\lambda_i+m\hbar)}\right).\]
We are assuming that there is at least one negative line bundle.
We expand the second factor of $I_V^{eq}$ as a polynomial of
$\hbar^{-1}$. Each negative line bundle produces a factor of
$\displaystyle{\frac{p}{\hbar}}$. For example, in the case
$V=\mathcal O(k)\oplus \mathcal O(-l)$ the expansion yields:
\begin{equation}
I_V^{eq}=\text{exp}\left(\frac{t_0+pt}{\hbar}\right)\left(1+\frac{p}{\hbar}\sum_{d=1}^{\infty}q^d\frac{(-1)^{ld}(ld-1)!(kd)!}{(d!)^{s+1}}\frac{1}{\hbar^{d(s+1-k-l)}}+o(\frac{1}{\hbar^2})\right)
\end{equation}
If $V$ contains two or more negative line bundles it follows that
\[\displaystyle{I_V^{eq}=\text{exp}\left(\frac{t_0+pt}{\hbar}\right)\left(1+o(\frac{1}{\hbar^2})\right)}\]
Lemma $(5.4.1)$ and $(5.4.2)$ imply that $J_V^{eq}=I_V^{eq}$. If
$\sum_{i\in I}k_i+\sum_{j\in J}l_j<s+1$ the presence of
$\displaystyle{\frac{1}{\hbar^{d(s+1-k-l)}}}$ in the above
expansion shows that again
\[\displaystyle{I_V^{eq}=\text{exp}\left(\frac{t_0+pt}{\hbar}\right)\left(1+o(\frac{1}{\hbar^2})\right)}\] hence $J^{eq}=I^{eq}$.
We may assume that $\sum_{i\in I}k_i+\sum_{j\in J}l_j=s+1$ and
$|J|=1$. In this case
\[I_V^{eq}=\text{exp}\left(\frac{t_0+pt}{\hbar}\right)\left(1+I_1\frac{p}{\hbar}+o(\frac{1}{\hbar^2})\right)\] where
$I_1$ is a power series of $q$ whose first term is zero. For
example if $V=\mathcal O(k)\oplus \mathcal O(-l)$ the power series
$I_1$ is
\[I_1=\sum_{d=1}^{\infty}q^d\frac{(-1)^{ld}(ld-1)!(kd)!}{(d!)^{s+1}}\]
Recall that $S=1+o(\hbar^{-2})$. Therefore
\begin{equation}
\text{exp}(\frac{I_1p}{\hbar})S(qe^{I_1},\hbar)=1+I_1\displaystyle
\frac{p}{\hbar}+o(\hbar^{-2}).
\end{equation}
Lemma $5.4.2$ implies that both $\displaystyle
\text{exp}\left(\frac{I_1p}{\hbar}\right)S(qe^{I_1},\hbar)$ and
$S'(q,\hbar)$ satisfy the conditions of the Lemma $5.4.1$. It
follows that
\[\displaystyle \text{exp}\left(\frac{I_1p}{\hbar}\right)S(qe^{I_1},\hbar)=S'(q,\hbar)\]
Multiplying both sides of this identity by $\displaystyle
\text{exp}\left(\frac{t_0+pt}{\hbar}\right)$ yields
\begin{equation}
J_V^{eq}(t_0,t+I_1)=I_V^{eq}(t_0,t).
\end{equation}
This completes the proof.$\dagger$

\vspace{0.2in}

\begin{corollary}
Let $V=(\oplus_{i\in I}{\mathcal O}(k_i))\oplus(\oplus_{j\in J}{\mathcal O}(-l_j))$. For $|J|>1$ or $k+l<s+1$
\[{e_1}_*\left(\frac{\EE'_d\EE^-_d}{\hbar(\hbar-c)}\right)=\frac{\prod_{i\in I}\prod_{m=1}^{k_id}(k_iH+m\hbar)\prod_{j\in J}\prod_{m=1}^{l_jd-1}(-l_jH-m\hbar)}{\prod_{m=1}^{d}(H+m\hbar)^{s+1}}.\]
\end{corollary}
{\bf Proof.} As mentioned above in this case we have $J_V^{eq}=I_V^{eq}$. Recall that
\[J_V^{eq}=\text{exp}\left(\frac{t_0+pt_1}{\hbar}\right)\left(1+\sum_{d>0}q^d{e_1}_*\left(\frac{E'_dE_d^-}{\hbar(\hbar-c)}\right)\cup \prod_{j\in J}(-l_jp)\right).\] We obtain the equivariant identity:
\[{e_1}_*\left(\frac{E'_dE_d^-}{\hbar(\hbar-c)}\right)\cup \prod_{j\in J}(-l_jp)=\frac{\prod_{i\in I}\prod_{m=1}^{k_id}(k_ip+m\hbar)\prod_{j\in J}\prod_{m=0}^{l_jd-1}(-l_jp-m\hbar)}{\prod_{m=1}^{d}\prod_{k=0}^{s}(p-\lambda_k+m\hbar)}.\] The restriction of $p$ to any fixed point $p_i$ is nonzero. This implies that $p$ is invertible. Therefore we obtain
\[{e_1}_*\left(\frac{E'_dE_d^-}{\hbar(\hbar-c)}\right)=\frac{\prod_{i\in I}\prod_{m=1}^{k_id}(k_ip+m\hbar)\prod_{j\in J}\prod_{m=1}^{l_jd-1}(-l_jp-m\hbar)}{\prod_{m=1}^{d}\prod_{k=1}^{s+1}(p-\lambda_k+m\hbar)}.\] The nonequivariant limit of this identity reads
\[{e_1}_*\left(\frac{\EE'_d\EE^-_d}{\hbar(\hbar-c)}\right)=\frac{\prod_{i\in I}\prod_{m=1}^{k_id}(k_iH+m\hbar)\prod_{j\in J}\prod_{m=1}^{l_jd-1}(-l_jH-m\hbar)}{\prod_{m=1}^{d}(H+m\hbar)^{s+1}}.\] The theorem is proven.$\dagger$

\smallskip

This corollary is particularly useful when $\text{Euler}(V^-)=0$
in $\PP^s$. In this case
$\displaystyle{J_V=I_V=\text{exp}\left(\frac{t_0+Ht_1}{\hbar}\right)}$
hence the mirror theorem is true trivially. An example of such a
situation is $V={\mathcal O}_{\PP^1}(-1)\oplus {\mathcal
O}_{\PP^1}(-1)$ which is treated in the next section.

\section{Examples}

\subsection{\bf Multiple covers} Let $C$ be a smooth rational curve in a Calabi-Yau threefold $X$ with normal bundle $N={\mathcal O}(-1)\oplus {\mathcal O}(-1)$ and
$\beta=[C]\in H_2(X,\ZZ)$. Since $K_X={\mathcal O}_X$ the expected
dimension of the moduli space ${\overline M}_{0,0} (X,d\beta)$ is
zero. However, this moduli space contains a component of positive
dimension, namely ${\overline M}_{0,0} (\PP^1,d)$. Indeed, let
$f:\PP^1\rightarrow C$ be an isomorphism, and $g:\PP^1\rightarrow
\PP^1$ a degree $d$ multiple cover. Then $f\circ g$ is a stable
map that belongs to ${\overline M}_{0,0} (X,d\beta)$. For a proof
of the fact that ${\overline M}_{0,0} (\PP^1,d)$ is a component of
${\overline M}_{0,0} (X,d\beta)$ see section $7.4.4$ in
\cite{[9]}. Let $N_d$ be the degree of $[{\overline M}_{0,0}
(X,d\beta)]^{\text{virt}}$. We want to compute the contribution
$n_d$ of ${\overline M}_{0,0} (\PP^1,d)$ to $N_d$. Kontsevich
asserted in \cite{[24]} and Behrend proved in \cite{[2]} that the
restriction of $[{\overline M}_{0,0} (X,d\beta)]^{\text{virt}}$ to
${\overline M}_{0,0} (\PP^1,d)$ is precisely $\EE_d$ for
$V={\mathcal O}(-1)\oplus {\mathcal O}(-1)$. Therefore:
\[n_d=\int_{{\overline M}_{0,0} (\PP^1,d)}\EE_d.\] Note that dim ${\overline M}_{0,0} (\PP^1,d)=2d-2$ and the
rank of the bundle $V_d$ is also $2d-2$. We use the mirror theorem
to compute numbers $n_d$. Since $V$ contains two negative line
bundles we can apply Corollary $5.5.1$
\[{e_1}_*\left(\frac{\EE_d}{\hbar(\hbar-c)}\right)=\frac{\prod_{m=1}^{d-1}(-H-m\hbar)^2}{\prod_{m=1}^{d}(H+m\hbar)^2}=\frac{1}{(H+d\hbar)^2}.\] An expansion of the left hand side using the divisor property for the modified gravitational descendants yields
\[{e_1}_*\left(\frac{\EE_d}{\hbar(\hbar-c)}\right)=\frac{dn_d}{\hbar^2}+\frac{H}{\hbar^3}\int_{{\overline M}_{0,1} (\PP^1,d)}c\cup \EE_d.\]
where $c$ is the chern class of the cotangent line bundle at the marked point. On the other hand:
\[\frac{1}{(H+d\hbar)^2}=\frac{1}{d^2\hbar^2}-\frac{2H}{d^3\hbar^3}.\]
We obtain the Aspinwall-Morrison formula
\[n_d=\frac{1}{d^3},\] which has been proved by several different methods \cite{[24]},{[28]},{[32]}.
We also obtain
\[\int_{{\overline M}_{0,1} (\PP^1,d)}c\cup \EE_d=-\frac{2}{d^3}.\]

\subsection{Virtual numbers of plane curves}
Let $X$ be a Calabi-Yau threefold containing a $\PP^2$. As we saw
in Remark $4.2.1$, the normal bundle of ${\PP}^2$ in $X$ is
$K_{{\PP}^2}={\mathcal O}(-3)$. Let $C$ be a rational curve of
degree $d$ in ${\PP}^2$. Since $K_X={\mathcal O}_X$, the expected
dimension of the moduli space ${\overline M}_{0,0} (X,[C])$ is
zero. Lemma~\ref{lemma: eqmod} says that ${\overline M}_{0,0}
({\PP}^2,d)={\overline M}_{0,0} (X,[C])$, hence the dimension of
this moduli stack is $3d-1$. Recall the diagram
\[ \begin{CD}
{\overline M}_{0,1} ({\PP}^2,d)@>e_1>>{\PP}^2 \\
@VV {\pi_1} V \\
{\overline M}_{0,0} ({\PP}^2,d)
\end{CD} \]
From Lemma~\ref{lemma: refeuler}, the virtual fundamental class of
${\overline M}_{0,0} (X,[C])$ is the refined top Chern class of
the bundle $V_d=R^1{\pi_1}_*(e^*_1(K_{{\PP}^2}))$ over ${\overline
M}_{0,0} ({\PP}^2,d)$. The zero pointed Gromov-Witten invariant:
\[N_d:=\text{deg}[{\overline M}_{0,0} (X,[C])]^{\text{virt}}=\int_{{\overline M}_{0,0} ({\PP}^2,d)} \EE_d\] is
called {\it the virtual number of degree $d$ rational curves in
$X$}. As promised in Remark $3.2.1$, we will show that the
modified equivariant quantum product in this case has a
nonequivariant limit. We will also use The Mirror Theorem to
calculate these numbers $N_d$.

The modified pairing on $\PP^2_{\CC^*}$ corresponding to
$V={\mathcal O}_{\PP^2}(-3)$ is
\[\langle a,b \rangle:=\int_{(\PP^2)_{\CC^*}}a\cup b\cup \left(\frac{1}{-3p-\lambda}\right).\]
Recall that $p$ denotes the equivariant hyperplane class in
$\PP^2$. Then $1,p,p^2$ is a basis for $\mathcal R$ as a
$\QQ(\lambda)$-module. A simple calculation shows that $-\lambda
p^2,-3p^2-\lambda H,-3p-\lambda$ is its dual basis with respect to
the above pairing. Since both bases and $E_d$ are polynomials in
$\lambda$, we can restrict $\tilde{I}_d$ and $*_V$ in ${\mathcal
P}=H^*(\PP^2,\QQ[\lambda])$ and take the nonequivariant limit.
Denote by $H$ the nonequivariant limit of $p$. We obtain the
following nonequivariant quantum product on
$H^*{\PP}^2\otimes{\QQ}[[q]]$
\[a*_Vb:=a\cup b+\sum_{d=1}^{\infty}q^dT^kI_d(a,b,-3HT_k)\] where $T^0=1, T^1=H , T^2=H^2$ and
for $\gamma_1,\gamma_2,...,\gamma_n\in H^*{\PP}^2$
\begin{equation}
I_d(\gamma_1,\gamma_2,...,\gamma_n)=\int_{{\overline M}_{0,n} ({\PP}^2,d)}e_1^*{\gamma_1}\cup e_2^*{\gamma_2}\cup...\cup e_n^*{\gamma_n}\cup \EE_d.
\end{equation}
For example, using the divisor axiom we obtain
\[H*_VH={H^2}(1-3\sum_{d>0}q^dd^3N_d).\]
Theorem~\ref{theorem: ring} implies the following:
\begin{theorem}
$(H^*{\PP^2}\otimes{\QQ}[[q]],*_V)$ is an associative, commutative
and unital ring with unity $1=[\PP^2]$.
\end{theorem}
Denote by $i$ the embedding $i:{\PP}^2\hookrightarrow X.$ Since
the normal bundle of $\PP^2$ in $X$ is ${\mathcal O}_{\PP^2}(-3)$,
it follows that $\displaystyle{i^*(-\frac{1}{3}[\PP^2])=T^1}$ and
$\displaystyle{i^*(-\frac{1}{3}[l])=T^2.}$ Therefore the map
$i^*:(H^*X,\QQ)\rightarrow (H^*{\PP}^2,\QQ)$ is surjective.
Consider the small quantum cohomology rings
$SQH^*X=(H^*X\otimes{\QQ}[[\beta]],*)$ and
$SQH^*_V\PP^2:=(H^*{\PP}^2\otimes{\QQ}[[q]],*_V)$ where the
products are given by three point correlators. Recall from section
$4.2$ the extension of $i^*:H^*(X,\QQ)\rightarrow H^*(\PP^2,\QQ)$
to $\tilde{i^*}:SQH^*X\rightarrow SQH^*_V\PP^2.$ There is a
natural relation between the modified quantum product in $\PP^2$
and the pure product in $X$.
\begin{theorem}
The map $\tilde{i^*}$ is a ring homomorphism.
\end{theorem}
{\bf Proof.} Complete $\tau^0=[X],\tau^1=-\frac{1}{3}[\PP^2],\tau^2=-\frac{1}{3}H$ into a basis of $(H^*X,\QQ)$ by adding elements from  $\text{Ker}(i^*)$.
Let $\tau_0=[pt],\tau_1=H,\tau_2=[\PP^2],...$ be the dual basis. Let $a,b\in H^*X$. We want to show
\[\tilde{i}^*(a*b)=i^*(a)*_Vi^*(b).\] But
\[a*b=\sum_{\beta\in H_2(X,\QQ) }\sum_{r}q^{\beta}\tau^r\int_{[{\overline M}_{0,3}(X,\beta)]^{\text{virt}}}e_1^*a\cup e^*_2b\cup e_3^*\tau_r.\] Note that this formula is true for a $\ZZ$-basis, but due to the uniqueness of the quantum product, it is true for any $\QQ$-basis as well. Now,
\[\tilde{i^*}(a*b)=\sum_{d\geq 0}\sum_{r}q^di^*\tau^r\int_{{\overline M}_{0,3}(\PP^2,d)}e_1^*(i^*a)\cup e^*_2(i^*b)\cup e_3^*(i^*\tau_r)\EE_d.\]
But $i^*(\tau^r)=T^k$ for $r=0,1,2$ and $i^*(\tau^r)=0$ for $r\geq 2$. The theorem follows from the readily checked fact: $i^*(\tau_k)=-3HT_k$ for $k=0,1,2$.$\dagger$

\vspace{0.2in}

Using the divisor and fundamental class properties of the modified
gravitational descendants it is easy to show that:
\[J_V=\text{exp}\left(\frac{t_0+t_1H}{\hbar}\right)\left(1-3\frac{H^2}{{\hbar}^2}\sum_{d=1}^{\infty}q^ddN_d\right).\]
The hypergeometric series corresponding to the total space of
$V={\mathcal O}_{\PP^2}(-3)$ is:
\[I_V:=\text{exp}\left(\frac{t_0+t_1H}{\hbar}\right)\sum_{d=0}^{\infty}q^d\frac{\prod_{m=0}^{3d-1}(-3H-m\hbar)}{\prod_{m=1}^{d}(H+m\hbar)^3}.\]
We expand this function
\[I_V=\text{exp}\left(\frac{t_0+t_1H}{\hbar}\right)\left(1+I_1\frac{H}{\hbar}+o(\frac{1}{\hbar})\right)\] where \[I_1=3\sum_{d=1}^{\infty}q^d(-1)^d\frac{(3d-1)!}{(d!)^3}.\]
The mirror theorem for this case says that $J(t_0,t_1+I_1)=I_V(t_0,t_1)$. This theorem allows us to compute the virtual number of rational plane curves in the Calabi-Yau $X$. The first few numbers are $3, \displaystyle\frac{-45}{8}, \displaystyle\frac{244}{9}$.

\vspace{+10 pt}
\noindent Department of Mathematics and Statistics
\noindent American University
\noindent 4400 Massachusetts Ave
\noindent Washington, Dc 20016
\noindent aelezi@american.edu

\end{document}